\documentclass{amsart}
\usepackage[english]{babel}
\usepackage{amsmath}
\usepackage[table]{xcolor}
\makeatletter
\renewcommand{\@eqnnum}{\hb@xt@.01\textwidth{\hfil\theequation}}
\makeatletter

\usepackage{graphicx}
\graphicspath{{Image/}}
\usepackage[foot]{amsaddr}
\usepackage{xcolor}
\usepackage{subcaption}
\usepackage{float}
\usepackage{multirow}
\newcommand{\fig}[1]{Fig.~\ref{#1}}{\color{blue}}

\usepackage{comment}
\usepackage{stackengine,graphicx}
\usepackage{tikz}
\usepackage[top=1in, bottom=1.25in, left=1.25in, right=1.25in]{geometry}
\usepackage{overpic}

\usepackage{hyperref}
    \hypersetup{
        colorlinks   = true,
        allcolors   =blue    
    }
\usepackage{tikz}
\usetikzlibrary{calc,fit,backgrounds}
\usetikzlibrary{shapes.geometric, arrows}
\usetikzlibrary{chains,
                fit,
                positioning,
                shapes}

\def\u{{\bm u}}

\def\0{\boldsymbol{0}}

\newtheorem{rem}{Remark}[section]

\newcommand{\bm}[1]{\mbox{\boldmath{$#1$}}}

\begin{document}

% \title[\anna{UPDATE} A comparison of data-driven ROMs for atmospheric flow]{ \anna{UPDATE} A comparison of data-driven Reduced Order Models for the simulation of mesoscale atmospheric flow}

\title[Combining E-CAEs and RC Computing for Accurate Predictions of Atmospheric Flows]{Combining Extended Convolutional Autoencoders and Reservoir Computing for Accurate Reduced-Order Predictions of Atmospheric Flows}

\author{ Arash Hajisharifi$^1$, Michele Girfoglio$^2$, Annalisa Quaini$^{3}$, and Gianluigi Rozza$^{1,*}$}
\address{$^1$ mathLab, Mathematics Area, SISSA, via Bonomea 265, I-34136 Trieste, Italy%; Arash.hajisharifi@sissa.it, mgirfogl@sissa.it, grozza@sissa.it 
}
\address{$^2$ Department of Engineering, University of Palermo, Viale delle Scienze, Ed. 7, 90128 Palermo, Italy}
\address{$^3$ Department of Mathematics, University of Houston, 3551 Cullen Blvd, Houston TX 77204, USA}
\address{$^*$ Corresponding author: \url{grozza@sissa.it}}

\begin{abstract}
Forecasting atmospheric flows with traditional discretization methods, also called full order methods (e.g., finite element methods or finite volume methods), is computationally expensive. We propose to reduce the computational cost with a Reduced Order Model (ROM) that combines Extended Convolutional Autoencoders (E-CAE) and Reservoir Computing (RC). Thanks to an extended network depth, the E-CAE encodes the high-resolution data coming from the full order method into a compact latent representation and can decode it back into high-resolution with $75 \%$ lower reconstruction error than standard CAEs. The compressed data are fed to an RC network, which predicts their evolution. The advantage of RC networks is a reduced computational cost in the training phase compared to conventional predictive models.

We assess our data-driven ROM through well-known 2D and 3D benchmarks for atmospheric flows. We show that our ROM accurately reconstructs and predicts the future system dynamics with errors below $6\%$ in 2D and $8\%$ in 3D, while significantly reducing the computational cost of a full order simulation.
Compared to other ROMs available in the literature, such as Dynamic Mode Decomposition and 
Proper Orthogonal Decomposition with Interpolation, our ROM is as efficient but more accurate. 
Thus, it is a promising alternative to high-dimensional atmospheric simulations.

% \anna{Dobbiamo essere piu' precisi nel quantificare questo:}
% % ability to accurately reconstruct and predict the future system dynamics while maintaining the low computational cost. 
% \anna{se lo lasciamo cosi' non e' molto informativo e magari uno decide di non leggere l'articolo perche' l'abstract non e' convincente. Quantifichiamo cosa si intende per accurate e per efficient in 2D e 3D.}

% Compared to existing ROMs, our method provides an accurate and efficient solution for this problem.
%\anna{Anche quest'ultima frase va precisata, al momento e' troppo generica.}

% Compared to existing ROMs, our method achieves a higher compression ratio, improved forecasting accuracy, and faster inference speed, making it an efficient and accurate solution for nonlinear reduced-order modeling

% These results establish E-CAE + RC as an efficient and accurate alternative for nonlinear reduced-order modeling and atmospheric flow prediction, paving the way for future advancements in data-driven ROMs

\end{abstract}

\maketitle


\section{Introduction}

Traditionally, the accurate forecast of atmospheric flow dynamics has required multiple large-scale, time-dependent simulations based on classical, high-fidelity discretization methods
(e.g., finite element methods or finite volume methods). Despite a continuous increase in  computational power and tremendous advancements in such discretization methods, also called Full Order Models (FOMs), these simulations remain 
computationally expensive. The search for accurate
alternative approaches with a reduced computational has motivated a large body of literature in recent years.
For example, methods borrowed from 
Machine Learning have been applied to global weather forecasting. See, e.g., \cite{pathak2018model,rasp2021data,  schultz2021can,weyn2019can,GraphCast,FourCastNet,Pangu-Weather}. 

We tackle this challenge starting from a reduced order modeling perspective. Reduced Order Models (ROMs) have emerged as a powerful approach to reduce the computational cost of FOM simulations. In ROMs, the reduction in the cost is achieved 
by approximating the high-dimensional FOM solution in a lower dimensional space, while maintaining the essential feature of the system in its compressed representation.  This makes ROMs advantageous for exploring a wide range of physical parameters or conducting long-term forecasts. 
A ROM model is constructed in two phases. In the first phase, called {offline}, a comprehensive dataset of FOM solutions, corresponding to specific time instances or/and sample physical parameters, is collected. 
This dataset is used to generate a reduced basis containing the compact representation of the system dynamics. The offline phase is computationally expensive, but it is performed only once.
In the second phase, called {online}, the offline-generated reduced basis is used to quickly compute the solution for a new time instance or parameter value. This is an efficient approach for multi-query contexts arising from needs such as assessing the uncertainty in the computed solution or solving an inverse problem for parameter identification or optimization. For further details, the
readers is referred to, e.g., \cite{peter2021modelvol1,benner2020modelvol2,benner2020modelvol3,hesthaven2016certified, malik2017reduced,rozza2008reduced, benner2017model, benner2015survey, siena2024accuracy}.

ROMs can be classified into intrusive and non-intrusive strategies \cite{rozza2022advanced, xiao2019domain}. The intrusive ROMs are physics-based approaches where the governing equations of the high-fidelity system are projected onto the reduced space spanned by the reduced basis computed in the offline phase. For these methods, which are also
referred to as projection-based methods, one needs access to
the FOM solver to perform the projection onto the
reduced space \cite{carlberg2017galerkin, reyes2020projection, rowley2004model, swischuk2019projection, benner2015survey}, hence the name intrusive. On the other hand, non-intrusive ROMs rely only on the FOM solution data and for this reason they are also called data-driven. In other words, they learn the reduced space and approximate the system's dynamics from the observed data \cite{schlegel2015long, xiao2014non, osth2014need}.

In this paper, we propose a non-intrusive ROM. We have preferred a non-intrusive approach because intrusive ROMs for highly nonlinear problems require hyper-reduction 
techniques, which are both computationally expensive and problem-specific \cite{barrault2004empirical, chaturantabut2010nonlinear}. For such problems, 
non-intrusive ROMs tend to be more computationally efficient.
Additionally, non-intrusive ROMs can leverage advanced machine learning techniques, e.g., deep learning-based approaches. These approaches have been successfully used in, e.g., FourCastNet \cite{FourCastNet}, GraphCast \cite{GraphCast}, and Pangu-Weather \cite{Pangu-Weather} to forecast complex atmospheric systems.

This work can be viewed as an extension of our previous paper \cite{hajisharifi2024comparison}, which compared three ``off-the-shelf'' non-intrusive ROMs, namely Dynamic Mode Decomposition (DMD) \cite{kutz2016dynamic,schmid2010dynamic, schmid2011applications, tu2013dynamic}, 
Hankel Dynamic Mode Decomposition (HDMD)
\cite{arbabi2017ergodic,  curtis2023machine, fujii2019data, jiang2015study}, and 
Proper Orthogonal Decomposition with Interpolation (PODI) \cite{shah2022finite, hajisharifi2023non, HAJISHARIFI2024106361,demo1,ripepi2018reduced}.
In \cite{hajisharifi2024comparison}, we showed 
that, although DMD and HDMD are designed to predict the dynamics of a system, their accuracy deteriorates 
rather quickly as the forecast time window expands.
On the other hand, the PODI solution is accurate for 
the entire duration of the time interval of interest thanks to the use of interpolation with radial basis functions. However, the interpolatory nature of PODI
limits its usefulness for predictions. In general, 
the strong limitation of all three methods is that they
use \emph{linear maps} to represent the high-dimensional FOM solution in a lower-dimensional spaces.
As a result, DMD, HDMD, and PODI often struggle to reproduce and predict strongly non-linear dynamics.

%In the geophysical fluid dynamics community, POD, often referred to as Empirical Orthogonal Function (EOF) analysis, is widely used to identify spatio-temporal coherent meteorological structure such as the Madden-Julian Oscillation, the Quasi-Biennial Oscillation and the El Niño-Southern Oscillation \cite{lario2022neural, pawar2022equation, schmidt2019spectral} which are characterized by large spatial and temporal scales. Although EOF analysis has been demonstrated to effectively identify the time evolution of the global scale, it is mostly constrained to the system identification. Recently, \cite{pawar2022equation} proposed data-driven ROMs based on EOF analysis to forecast the weekly average sea surface temperature. Hajisharifi et al. \cite{hajisharifi2024comparison} argued that PODI performs well in the interpolatory regimes for a parametric study of atmospheric flow problems and it can accurately identifies the underlying dynamics. 
%optimized DMD (ODMD) \cite{askham2018variable, sashidhar2022bagging} 

Deep-learning-based ROM approaches overcome
this limitation by adopting \emph{nonlinear maps}
to compress the representation of nonlinear dynamics. 
Several studies (e.g., \cite{lee2020model, fresca2021comprehensive, halder2022non, maulik2021reduced}) have shown that these ROMs outperform existing linear ROMs in terms of accuracy. 
Autoencoders (AE) are a type of neural network 
widely used  to reduce the dimensionality of 
high-dimensional data, typically from images, 
by learning nonlinear representations  \cite{hinton2006reducing, eivazi2020deep}. 
In our case, the images are the plots of FOM solutions. 
The  AE learns to encode the FOM solutions 
into lower-dimensional latent space representation by extracting the most essential features in the data. Then, a decoder can be applied to reconstruct the data back in the original high-dimensional space from the compressed representation. 
Convolutional Autoencoders (CAE) \cite{mao2016image, zhang2018better, fries2022lasdi} replace fully connected layers in standard AE with convolutional layers, which can capture finer details and more features of spatial correlations and thus improve the quality of compressed 
data representation \cite{halder2022non, gonzalez2018deep, omata2019novel}.

We propose a novel nonlinear ROM approach that integrates an enhanced version of CAE to efficiently capture the spatial correlations (i.e., higher compression ratio and lower reconstruction error than a standard CAE), with a Reservoir Computing (RC) framework to learn the evolution of the CAE latent spaces and accurately predict their future dynamics.

RC networks, derived from Recurrent Neural Network (RNN) theory \cite{jaeger2004harnessing}, process sequential data. Unlike the CAE that compresses 
the high-resolution FOM data (dimension of the order of tens of thousands or millions) into a compact latent space (dimension of order ten) and reconstructs it back to the original 
resolution, the RC does not reconstruct the high-resolution data. Instead, it maps the input data, which in our case is the latent space representation given by the CAE, to a higher dimensional reservoir state space (dimension of the order of hundreds) through reservoir nodes, with the goal of better capturing the temporal dependency in the data.  
% The RC framework maps the input data, which in our case is the latent space representation given by the CAE, to a high-dimensional state space thrgouh reservoir nodes, with the goal of better capturing the temporal dependency in the data. 
%\anna{(secondo me qui dobbiamo dare ordini di grandezza perche' c'e' l'high dimensiona data che viene dato al CAE e quello che gestisce l'RC e non sono la stessa dimensione)} 
Unlike other RNNs, the weights of the recurrent connections are not trained but initialized randomly and remain fixed throughout the process. The only trainable parameters in RC networks are the output layer weights.
These are trained using a simple linear regression algorithm. This simple  training process significantly reduced the computational cost of RC compared to other predictive  networks \cite{schrauwen2007overview, lukovsevivcius2009reservoir, tanaka2019recent}, like Long-Short Time Memory (LSTM) \cite{hochreiter1997long,maulik2021reduced, pawar2020data, li2019deep, HAJISHARIFI2024106361},  Gated Recurrent Unit (GRU)
\cite{cho2014learning}, and Transformers
\cite{vaswani2017attention,geneva2022transformers, hassanian2023deciphering, sarma2024prediction, yousif2023transformer}.
In fact, LSTMs and GRUs are usually slow to train because they are sequential and, consequently, cannot exploit parallel architectures.
Transformers employ an attention mechanism to avoid the sequential data processing. This mechanism enables them to effectively model long-term dependencies of a dynamical system. Despite their ability to take advantage of parallel computing, 
transformers are computationally expensive and memory-intensive, which makes them less efficient for real time applications.
Through the combination of CAE and RC, our ROM aims to improve prediction accuracy, while keeping the computational cost low.

An alternative to nonlinear ROMs, such as the one presented in this work, is to improve the 
representation of nonlinear dynamics by a linear 
ROM though data augmentation. Data augmentation can be done using, e.g., techniques borrowed from 
optimal transport theory \cite{khamlich2024optimal}. 

%LSTM is introduced to improve the performance of RNNs by solving the vanishing gradient problem. Similar to LSTMs,  GRUs tackle the vanishing gradient issue by introducing a different strategy called the Gated Recurrent Unit . Despite the advantages, these networks are usually slow to train as they are sequential and consequently, they can not exploit parallel architectures. 
%To tackle this issue, transformer networks were introduced. Transformers employ an attention mechanism to avoid the sequential data processing that enables them to effectively model long-term dependencies of a dynamical system. Despite their ability to take advantage of parallel computing, Transformers have been shown to outperform other prediction methods \cite{yousif2023transformer, borrelli2022predicting}.  However, they are computationally expensive and memory-intensive that makes them less efficient for real time applications. Moreover, they require large dataset for training as they tend to easily overfit on smaller datasets \cite{sanford2024representational, peng2024limitations}.

%\textcolor{blue}{Non e' meglio di spiegare anche un po altre reti predittive ?}

%\anna{Arash si occpua di controllare che citiamo tutti i lavori citati anche da Moaad, citare il lavoro con Mooad e quello di Peherstorfer, controllare se gli item della biblio in arXiv sono apparsi in giornali. }
%\textcolor{blue}{Ho aggiunto il lavoro di Mooad, anche altre three ref che era correlato al nostro lavoro}

The rest of the paper is organized as follows. In  Sec.~\ref{sec:FOM}, we state the compressible Euler equations for low Mach, stratified flows, which represent our FOM.  Sec.~\ref{sec:ROM} describes the details of the proposed ROM approach. 
Sec.~\ref{sec:res} reports the results for three well-known benchmarks: 2D and 3D rising thermal bubble and 2D density current. Finally, Sec.~\ref{sec:conc} provides conclusions and future perspectives.

\section{The full order model}\label{sec:FOM}

We are interested in the dynamic of the dry atmosphere (i.e.,
no humidity) and the effects of solar radiation and ground heat flux are neglected for simplicity. In a spatial domain of interest $\Omega$ over the time interval $(0, t_f]$, the Euler equations for dry air dynamics read: find air density $\rho$, wind velocity $\u$,
and total energy $e$ such that
\begin{align}
%\centering
&\frac{\partial \rho}{\partial t} + \nabla \cdot ( \rho \u ) = 0 &&\text{in} \ \Omega \ \times \ (0, t_f], \label{continuity} \\
&\frac{\partial (\rho \u)}{\partial t} + \nabla \cdot ( \rho \u \otimes \u) + \nabla p + \rho g \widehat{\mathbf k} = \boldsymbol{0}  &&\text{in} \ \Omega \ \times  \ (0, t_f], \label{momentu} \\
&\frac{\partial (\rho e)}{\partial t} + \nabla \cdot ( \rho \u e) + \nabla \cdot (p \u) = 0  &&\text{in} \ \Omega \ \times  \ (0, t_f],
\label{energy}
\end{align}
where $p$ is the pressure, $\widehat{\mathbf k}$ is the unit vector pointing the vertical axis $z$, and $g$ is the gravitational constant. We note that eq.~\eqref{continuity}, 
\eqref{momentu}, and \eqref{energy} describe conservation
of mass, momentum and energy, respectively.  
The total energy density can be written as follows:  
\begin{align}
e = c_v T + K + g z,\quad K = |\u|^2/2
\label{eq_total_energy}
\end{align}
where $c_{v}$ is the specific heat capacity at constant volume and $T$ is the absolute temperature. 
Let $c_{p}$ is the specific heat capacity at constant pressure.
By introducing the specific enthalpy, $h = c_v T + p/ \rho  = c_p T$, the total energy density can be rewritten as:
\begin{align}
e = h - {p}/{\rho} + K + gz,
\label{eq_total_energy2}
\end{align}
where $K $ is the kinetic energy density defined in \eqref{eq_total_energy}. Then, by plugging
\eqref{eq_total_energy2} into \eqref{energy}
and accounting for \eqref{continuity},
we obtain:
\begin{align}
\frac{\partial (\rho h)}{\partial t} + \nabla \cdot ( \rho \u h) + \frac{\partial (\rho K)}{\partial t} + \nabla \cdot (\rho \u K) -\frac{\partial p}{\partial t} + \rho g \u \cdot \widehat{\mathbf k}= 0  \hspace{1cm} \text{in} \ \Omega \ \times  \ (0, t_f].
\label{energy2}
\end{align}

To close system \eqref{continuity}-\eqref{energy},
we need an equation of state. We assume that the dry atmosphere
behaves like an ideal gas, thus we have:
\begin{align}
p = \rho R T,
\label{EQ_of_State}
\end{align}
where $R$ is the specific gas constant of dry air.  

Let us to write the pressure as the sum of a fluctuation $p'$ with respect to a hydrostatic term:
\begin{equation}
    p=p_g + \rho g z + p', \label{EQ_p_splitting}
\end{equation}
where $p_g = 10^5$ Pa is the atmospheric pressure at the ground. By plugging \eqref{EQ_p_splitting} into \eqref{momentu}, the momentum conservation equation can be rewritten as:
\begin{align}
\frac{\partial (\rho \u)}{\partial t} + \nabla \cdot ( \rho \u \otimes \u) + \nabla p^\prime +  g z \nabla \rho = 0  \hspace{3cm} \text{in} \ \Omega \ \times  \ (0, t_f].
\label{momentu2}
\end{align}

We focus on formulation \eqref{continuity},\eqref{energy2},\eqref{EQ_of_State}-\eqref{momentu2} of the Euler equations.
Note that this system does not include a dissipation mechanism.
Hence, numerical errors are prone to get 
amplified over time, leading to numerical instabilities 
and eventually a simulation breakdown. There are several ways 
to address this problem. The easiest (and crudest) is
introduce an artificial dissipation term in \eqref{momentu2} and \eqref{energy2} with a constant artificial viscosity $\mu_a$:
\begin{align}
&\frac{\partial (\rho \u)}{\partial t} +  \nabla \cdot (\rho \u \otimes \u) + \nabla p' + gz \nabla \rho -  \nabla \cdot (2 \mu_a \boldsymbol{\epsilon}(\u)) + \nabla \left(\frac{2}{3}\mu_a \nabla \cdot \u \right)= 0 &&\text{in } \Omega \times (0,t_f],  \label{eq:mom_LES} \\
&\frac{\partial (\rho h)}{\partial t} +  \nabla \cdot (\rho \u h) + 
\frac{\partial (\rho K)}{\partial t} +  \nabla \cdot (\rho \u K) - \dfrac{\partial p}{\partial t}  +  
\rho g \u \cdot \widehat{\mathbf k}  - \nabla \cdot \left(\frac{\mu_a}{Pr} \nabla h \right) = 0 &&\text{in } \Omega \times (0,t_f],
\label{eq:ent_LES}
\end{align}
where $\boldsymbol{\epsilon}(\u) = (\nabla \u + (\nabla \u)^T)/2$ is the strain-rate tensor
and $Pr$ is the Prandtl number. In more sophisticated
approaches (see, e.g., \cite{clinco2023filter, girfoglio2023validationAIP, marras2015stabilized}), $\mu_a$ varies in space and time. To the effect of the method proposed
in this paper, it does not matter how $\mu_a$ is defined. Hence, 
we will keep it constant for simplicity. 
So, the starting point of our full order method is model 
\eqref{continuity},\eqref{EQ_of_State},\eqref{EQ_p_splitting},\eqref{eq:mom_LES},\eqref{eq:ent_LES}, plus suitable boundary conditions that will be specified for each test
problem considered in Sec.~\ref{sec:res}.

A quantity of interest for atmospheric studies is the potential temperature
\begin{align}
\theta = \frac{T}{\pi},  \quad \pi = \left(\frac{p}{p_g}\right)^{\frac{R}{c_p}},
\label{pot_temp}
\end{align}
i.e., the temperature that a parcel of dry air would have if it were expanded or compressed
adiabatically to standard pressure $p_g = 10^5$ Pa. % which is the atmospheric pressure at the ground. 
In \eqref{pot_temp}, $\pi$ is the so-called Exner pressure.
Similar to the splitting adopted for the pressure in \eqref{EQ_p_splitting}, we split the potential temperature 
into a hydrostatic value $\theta_0$, that depends only on the vertical coordinate $z$, and fluctuation $\theta^ \prime$ over it:
\begin{align}
\theta ^ \prime (x,y,z,t) = \theta (x,y,z,t) - \theta_0(z).
\label{pot_temp2}
\end{align}

To discretize
problem \eqref{continuity},\eqref{EQ_of_State},\eqref{EQ_p_splitting},\eqref{eq:mom_LES},\eqref{eq:ent_LES} in time and space, 
we proceed as follows. We choose a time step $\Delta t \in \mathbb{R}^+$ to divide time interval $(0,t_f]$ to $t^n = t_0 + n \Delta t$, with $n = 0, ..., N_{tf}$ and $t_f = 0 + N_{tf} \Delta t$. The time derivatives in \eqref{continuity}, \eqref{eq:mom_LES}, and \eqref{eq:ent_LES} are discretized adopting the Euler scheme. The convective terms in eq.~\eqref{eq:mom_LES} and \eqref{eq:ent_LES} are handled semi-implicitly, while the diffusive terms are treated implicitly. On the other hand, eq.~\eqref{continuity} is solved explicitly.
For the space discretization, we adopt a second-order finite volume scheme. For this, the computational domain $\Omega$ is divided into cells or control volumes $\Omega_i$, with $i = 1, \dots, N_{c}$, where $N_{c}$ is the total number of cells in the mesh. 
To decouple the computation of the pressure from velocity, 
we adopt the PISO algorithm \cite{patankar1983calculation,issa1986solution,moukalled2016finite}. 
All of the above choices are so that we can design an efficient splitting method for the complex coupled problem at hand. The reader interested in more details about our numerical approach is referred to \cite{girfoglio2023validationAIP,GIRFOGLIO2025106510}.

The above full order model is implemented in GEA (Geophysical and Environmental Applications) \cite{GEA,GQR_GEA, clinco2023filter,girfoglio2023validationAIP}, an open-source package for atmosphere and ocean modeling based on the finite volume C++ library OpenFOAM\textsuperscript{\textregistered}.

It is important to stress that, although we 
have made specific choices for the full order method, the reduced order model
presented in the next section can be applied to
data generated by any other full order method.

%%%%%%%%%%%%%%%%%%%%%%%%%%%%%%%%%%%%%%%%%%%%%%%%%%%%%%%%%%%%%%%%%%%%%%%%%%%%%%%%%%%%%%%%
%%%%%%%%%%%%%%%%%%%%%%%%%%%%%%%%%%%%%% 

\section{The reduced order model}\label{sec:ROM}

% \begin{figure}[t]
%     \vspace{1cm}
%     %\begin{subfigure}{\linewidth} 
%     \centering
%     \begin{overpic}[width=0.95\textwidth, grid=false]{CNN.pdf}  
%          % \put(6.4,23){ROM}
%          % \put(-5,10){\tiny{\textcolor{black}{$128 \times 32$}}}
         
%          \put(4,2.5){\tiny{\textcolor{black}{$f_1$}}}
%          \put(17.5,5.5){\tiny{\textcolor{black}{$f_2$}}}
%          \put(29,9){\tiny{\textcolor{black}{$ f_3$}}}
%          \put(37,12){\tiny{\textcolor{black}{$ f_4$}}}

%          \put(51,27){\tiny{\textcolor{black}{flatten}}}
%          \put(51,25.5){\tiny{\textcolor{black}{layer}}}

%          \put(53.5,23){\tiny{\textcolor{black}{Latent}}}
%          \put(53.5,21.5){\tiny{\textcolor{black}{space}}}

%          \put(62,27){\tiny{\textcolor{black}{flatten}}}
%          \put(62,25.5){\tiny{\textcolor{black}{layer}}}

%          \put(62,12){\tiny{\textcolor{black}{$f_4$}}}
%          \put(68,9){\tiny{\textcolor{black}{$f_3$}}}  
%          \put(75,5.5){\tiny{\textcolor{black}{$f_2$}}}         
%          \put(87,2.5){\tiny{\textcolor{black}{$f_1$}}}
         
%       \end{overpic}
%       \captionsetup{list=no} 
%   \caption{Convolutional Autoencoder (CAE) architecture. The values $f_1$ to $f_4$  indicate the number of filters applied at each resolution level. }
%    \label{fig:CAE}
% \end{figure}

Despite being efficient, the full order method described in the previous section 
can be computationally expensive if one is interested in simulating atmospheric flows over large computational domains and/or for long periods of time. In fact, the FOM solution is
associated to a large number of degrees of freedom and, hence, it is high-dimensional.
Through the ROM,  we aim at approximating the FOM solution in a space with a reduced dimension while preserving the main dynamical features. 

We will present our method to obtain the ROM approximation for variable $\theta^\prime$:
\begin{align}
    \theta_h^\prime(\bm{x}, t) \approx \theta_r^\prime(\bm{x}, t),
\end{align}
where $\theta_r^\prime$ is the ROM approximation of the FOM solution $\theta_h^\prime$. However, the same method can be applied to find any other variable of interest. 

The preliminary step to find $\theta_r^\prime$ is to collect the so-called snapshots, i.e., the FOM solutions at a given time, $\theta_h^{\prime} (\bm{x},  t^i) \in \mathbb{R}^{N_c}$, with $i = 1, \dots, N_t$. The subsequent steps involve:
\begin{itemize}
    \item[-] a Convolutional Autoencoder (CAE), or its enhanced version, to encode the snapshots to a compressed representation, the so-called latent space representation;
    \item[-] a Reservoir Computing (RC) framework to model the temporal evolution of the latent space representation.
\end{itemize}
Then, the CAE decodes the predicted latent spaces to obtain the associated high-dimensional approximation, i.e., the actual quantity of interest.

In the subsections below, we provide more details about the CAE and the RC framework, which constitute the main building blocks of our ROM.

\subsection{Convolutional Autoencoders and Extended Convolutional Autoencoders}

Convolutional autoencoders \cite{halder2022non, gonzalez2018deep, maulik2021reduced} are neural networks that learn to map the input data to a lower representation \cite{hinton2006reducing, lusch2018deep}. 
They take in a two-dimensional spatial structured data instance, i.e., an image, which in our case plots a computed
solution, and process it until a one-dimensional vector representation  
is produced. This operation is performed by a convolutional neural network and it extracts meaningful features from the input image.
Since we are considering a time-dependent problem, our image evolves in time and so will its one-dimensional vector representation. 
It is logical to assume that another convolutional neural network could 
take in the one-dimensional vector representation and reconstruct
the image by, in some sense, reversing the feature extraction process.
A convolutional autoencoder performs both the feature extracting (encoder) 
and the image reconstruction (decoder) with convolutional neural networks. 
The convolutional layers enable the network to capture spatial patterns and extract correlations more effectively from high-dimensional data \cite{mao2016image, zhang2018better, gonzalez2018deep}.  %\anna{Possiamo dare una referenza per quest'ultima frase? Puo' essere anche generica, tipo
%un libro di testo. \textcolor{blue}{Done}}
%Between the encored and the decoder, there is an additional
%layer called bottleneck, which helps to compress the extracted features into a smaller vector representation.
% https://www.digitalocean.com/community/tutorials/convolutional-autoencoder

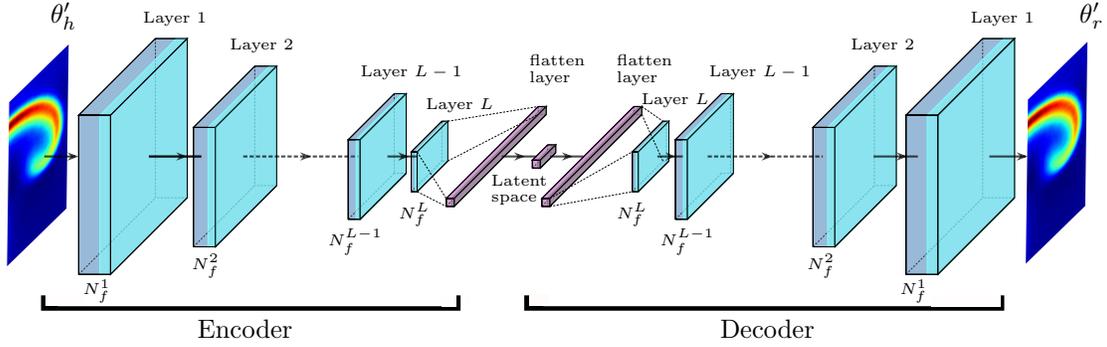
\begin{figure}[htb]
    \vspace{1cm}
    %\begin{subfigure}{\linewidth} 
    \centering
    \begin{overpic}[width=0.95\textwidth, grid=false]{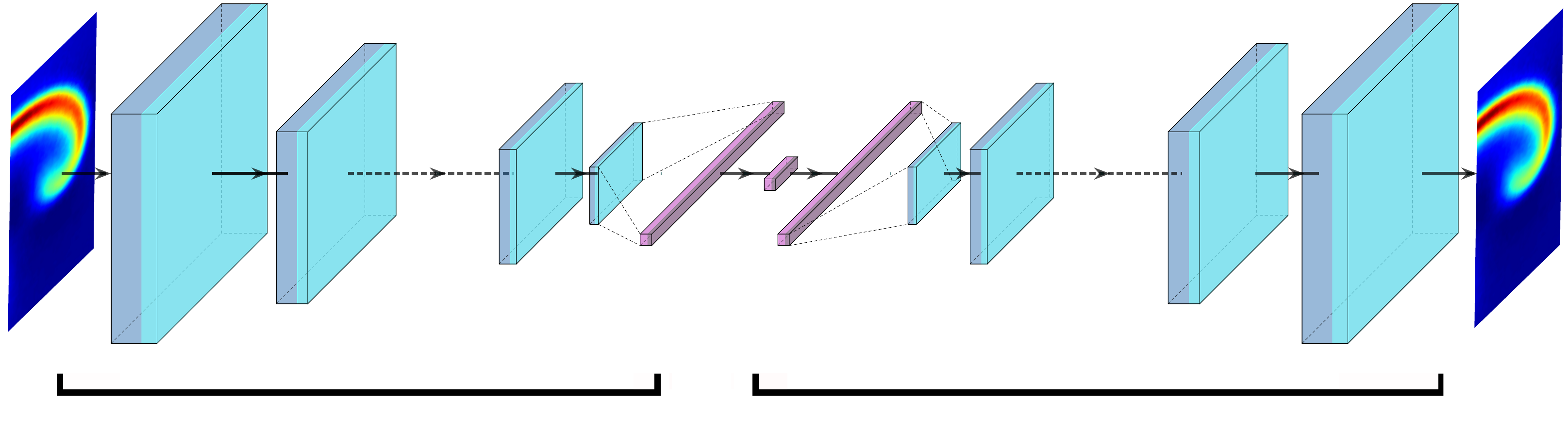}  
         % \put(6.4,23){ROM}
         % \put(-5,10){\tiny{\textcolor{black}{$128 \times 32$}}}
         \put(4.5,29){\textcolor{black}{$\theta ' _h$}}
         \put(7.5,4){\tiny{\textcolor{black}{$N_f^1$}}}
         \put(17.5,6.3){\tiny{\textcolor{black}{$N_f^2$}}}
         \put(30,8.5){\tiny{\textcolor{black}{$ N_f^{L-1}$}}}
         \put(36.5,11.){\tiny{\textcolor{black}{$ N_f^L$}}}

         \put(13,29){\tiny{\textcolor{black}{Layer 1}}}
         \put(21,26.5){\tiny{\textcolor{black}{Layer 2}}}
         \put(33,24){\tiny{\textcolor{black}{Layer $L-1$}}}
         \put(39,21){\tiny{\textcolor{black}{Layer $L$}}}

         \put(58.8,21.5){\tiny{\textcolor{black}{Layer $L$}}}
         \put(65,24){\tiny{\textcolor{black}{Layer $L-1$}}}
         \put(78,26.5){\tiny{\textcolor{black}{Layer 2}}}
         \put(89,29){\tiny{\textcolor{black}{Layer 1}}}         

         \put(48.5,25){\tiny{\textcolor{black}{flatten}}}
         \put(48.5,23.5){\tiny{\textcolor{black}{layer}}}

         \put(45,14){\tiny{\textcolor{black}{Latent}}}
         \put(45,12.5){\tiny{\textcolor{black}{space}}}

         \put(56.5,25){\tiny{\textcolor{black}{flatten}}}
         \put(56.5,23.5){\tiny{\textcolor{black}{layer}}}

         \put(56.5,11.){\tiny{\textcolor{black}{$N_f^L$}}}
         \put(60.5,8.5){\tiny{\textcolor{black}{$N_f^{L-1}$}}}  
         \put(74,6.3){\tiny{\textcolor{black}{$N_f^2$}}}         
         \put(82.5,4){\tiny{\textcolor{black}{$N_f^1$}}}
         \put(99,29){\textcolor{black}{$\theta ' _r$}}

         \put(18,0){\textcolor{black}{Encoder}}
         \put(66,0){\textcolor{black}{Decoder}}
         
      \end{overpic}
      \captionsetup{list=no} 
  \caption{CAE architecture with encorder and decoder that feature $L$ layers. For each layer $l$, with $l = 1, \dots, L$, 
  $N_f^l$ is the number of applied filters.}
   \label{fig:CAE}
\end{figure}

The encoder part of the CAE maps the input data $\theta_h^\prime(\bm{x}, t)$ 
to its latent space representation $\bm{z}(t) \in \mathbb{R}^{N_d}$, with $N_d \ll N_c$, through a sequence of convolutional layers $l =1, \dots, L$. See \fig{fig:CAE}. These layers gradually decrease spatial resolution to generate a compact representation in the final layer that captures the spatial pattern of the high-resolution input data. Formally, we can write the encoding process as:
\begin{align}
\bm{z}(t) = \mathcal{E}( \theta_h^\prime(\bm{x}, t)),
\label{equ_LS}
\end{align}
where $\mathcal{E}: \mathbb{R}^{N_c} \to \mathbb{R}^{N_d}$.
In each layer $l$, a set of filters is convolved over the input data to extract the spatial features. The number of filters applied in each layer is denoted with $N_f^l$ for $ l = 1, \dots, L$. The values of $L$ and $N_f^l$ 
are hyperparameters of the network that determine the feature extraction capacity overall and at each layer.

The decoder part of the CAE reconstructs the approximation 
$\theta_r^\prime(\mathbf{x}, t)$  from latent spaces $\bm{z}(t)$ using the transposed 
%\anna{(is it really simply the transposed? It's not some sort of approximation of the inverse?)}
convolution layers, which are called deconvolutional layers. See \fig{fig:CAE}. {This operation is not a direct mathematical inversion of the encoder. Instead, these layers learn a mapping to approximate the inverse process, i.e., to perform upsampling and generate the approximation of the FOM solution.} 
% i.e., to provide an approximation of the FOM solution,
Formally, this process can be written as:
\begin{align}
 \theta_r^\prime(\bm{x}, t) = \mathcal{D}(\bm{z}(t)),
\label{equ_dec}
\end{align}
where $\mathcal{D}: \mathbb{R}^{N_d} \to \mathbb{R}^{N_c}$. 
%\anna{se capisco giusto, qui $\mathcal{D}$ approssima $\mathcal{E}^{-1}$, no?} \textcolor{blue}{Si, esatto}
Through an optimization process, the CAE is trained to minimize the reconstruction error between the input data and the decoder's output, thereby ensuring that the latent representation is able to capture the most essential nonlinear spatial features of the system.
%\anna{Com'e' che si assicura questo?}.

Although standard CAEs are generally effective in capturing spatial patterns in image, for images stemming from 
the numerical approximation of complex flows 
their capability is not sufficient to reconstruct the FOM solution.

\begin{figure}[htb]
    \vspace{1cm}
    %\begin{subfigure}{\linewidth} 
    \centering
    \begin{overpic}[width=0.95\textwidth, grid=false]{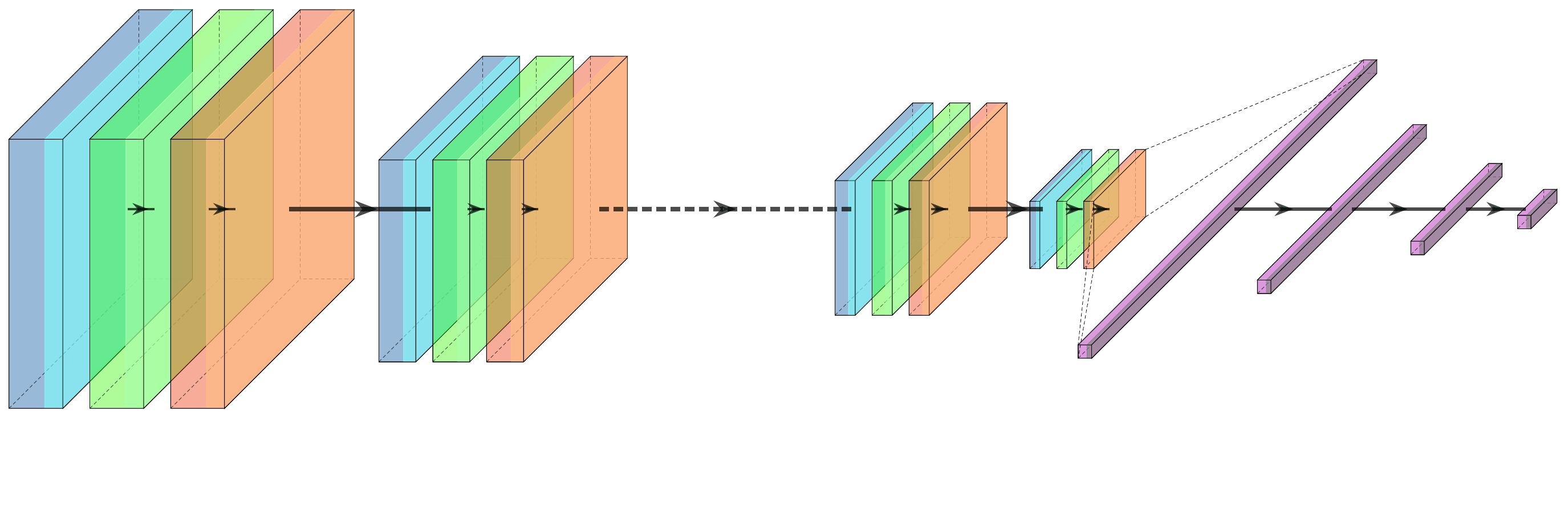}  
         % \put(6.4,23){ROM}
         % \put(-5,10){\tiny{\textcolor{black}{$128 \times 32$}}}
         
         \put(3,3.5){\tiny{\textcolor{black}{$n_f \times N_f^1$}}}
         \put(25,6){\tiny{\textcolor{black}{$n_f \times N_f^2$}}}
         \put(52.5,9){\tiny{\textcolor{black}{$n_f \times N_f^{L-1}$}}}
         \put(64,13){\tiny{\textcolor{black}{$n_f \times N_f^{L}$}}}

         \put(13,34){\tiny{\textcolor{black}{Layer 1}}}
         \put(32,30.5){\tiny{\textcolor{black}{Layer 2}}}
         \put(56,27){\tiny{\textcolor{black}{Layer $L-1$}}}
         \put(68,24){\tiny{\textcolor{black}{Layer $L$}}}

         \put(86,32){\tiny{\textcolor{black}{flatten}}}
         \put(86,30.5){\tiny{\textcolor{black}{layer}}}
         
         \put(89.5,27.5){\tiny{\textcolor{black}{Dense}}}
         \put(89.5,26){\tiny{\textcolor{black}{layer}}}

         \put(93.5,24.5){\tiny{\textcolor{black}{Dense}}}
         \put(93.5,23){\tiny{\textcolor{black}{layer}}}

         \put(99,23){\tiny{\textcolor{black}{Latent}}}
         \put(99,21.5){\tiny{\textcolor{black}{space}}}

      \end{overpic}
      \captionsetup{list=no} 
  \caption{Encoder architecture of the E-CAE with $L$ layers. For each layer $l$, with $l = 1, \dots, L$, we  apply $n_f$ sets of $N_f^l$ filters. %\anna{Figure to be updated \textcolor{blue}{done}}
  }
   \label{fig:E-CAE}
\end{figure}
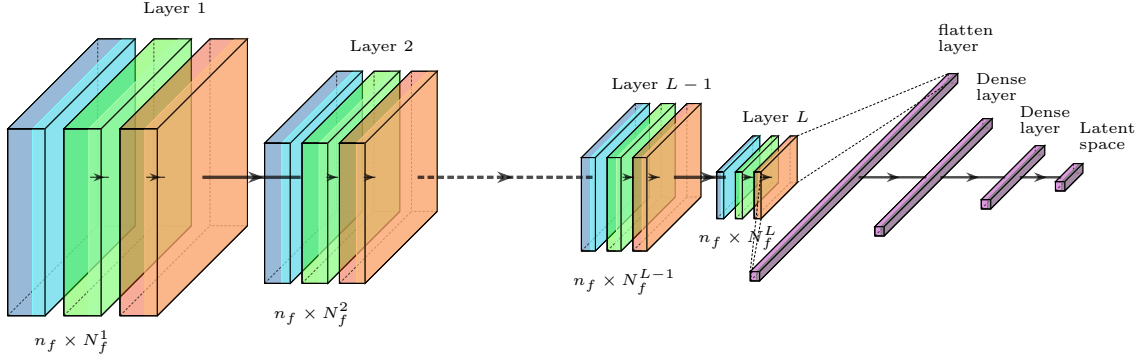

\label{sec:E-CAE}
To overcome this limitation, we consider a new architecture called {Extended Convolutional Autoencoder (E-CAE)}. 
In this network, which is inspired by modern deep learning architectures such as VGGNet\cite{simonyan2014very, swapna2020cnn}, ResNet \cite{he2016deep}, and DenseNet \cite{huang2017densely},
additional convolutional layers are stacked on top of the primary layer at each resolution level, i.e., at every layer $l$, with $l = 1, \dots, L$, we apply $n_f$ sets of $N_f^l$ filters. 
See \fig{fig:E-CAE}. This is meant to increase the network depth
and enhance feature extraction capability, so that
the network can learn more complex patterns before 
encoding them into the latent space.
For the results in Sec.~\ref{sec:res}, we will use the same value of $n_f$ for all the layers.
Because more features will be extracted by the E-CAE, 
we add so-called dense layers between the flatten layer
and the latent space to further compress the data and
make the overall dimensionality reduction smoother.
These dense layers are added %between the flatten layer
%and the latent space 
for both the encoder and 
the decoder and they are expected to help with faster convergence and lower reconstruction error. 

%This architecture has shown to be more capable of capturing more complex dynamics by improving feature extraction.

% \begin{figure}[htb!]
%     \centering \includegraphics[width=160mm,scale=0.6]{Depth_CNN_Arch.pdf}
%         \caption{CAE architecture}  
%         \label{fig:LS_RTB}
% \end{figure}

% \begin{figure}[htb!]
%     \vspace{1cm}
%     %\begin{subfigure}{\linewidth} 
%     \centering
%     \begin{overpic}[width=0.85\textwidth, grid=false]{Decoder.pdf}  
%          % \put(6.4,23){ROM}
%          % \put(-5,10){\tiny{\textcolor{black}{$128 \times 32$}}}
         
%          \put(2,0.5){\tiny{\textcolor{black}{$3 \times 256$}}}
%          \put(14,3){\tiny{\textcolor{black}{$3 \times 128$}}}
%          \put(28,5.5){\tiny{\textcolor{black}{$3 \times 64$}}}
%          \put(35.5,7){\tiny{\textcolor{black}{$3 \times 32$}}}

%          \put(57,7){\tiny{\textcolor{black}{$3 \times 32$}}}
%          \put(64,5.5){\tiny{\textcolor{black}{$3 \times 64$}}}  
%          \put(72.5,3){\tiny{\textcolor{black}{$3 \times 128$}}}         
%          \put(85,0.5){\tiny{\textcolor{black}{$3 \times 256$}}}
         
%       \end{overpic}
%       \captionsetup{list=no} 
%   \caption{E-CAE architecture.}
%    \label{fig:CAE2}
% \end{figure}

\subsection{Reservoir Computing}

\begin{figure}[htb]
    \vspace{1cm}
    %\begin{subfigure}{\linewidth} 
    \centering
    \begin{overpic}[width=0.6\textwidth, grid=false]{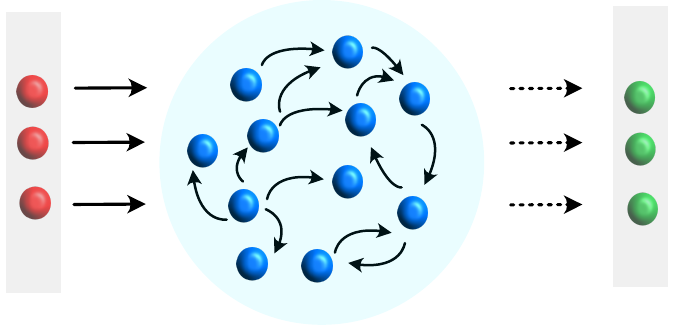}  
         % \put(6.4,23){ROM}
         % \put(-5,10){\tiny{\textcolor{black}{$128 \times 32$}}}
         
         \put(-2,50){\textcolor{black}{Input layer}}
         \put(35,50){\textcolor{black}{Reservoir node}}
         \put(84,50){\textcolor{black}{Output layer}}

         \put(12,10){\textcolor{black}{$\bm{W}_{in}$}}
         \put(45,2){\textcolor{black}{$\bm{W}_{res}$}}
         \put(76,10){\textcolor{black}{$\bm{W}_{out}$}}

      \end{overpic}
      \captionsetup{list=no} 
  \caption{Sketch of the three main components of an
  RC architecture: input layer, reservoir node, and output layer. With $\bm{W}_{in}$, $\bm{W}_{res}$ and $\bm{W}_{out}$, we denote the weight matrices that define the interactions among the components.}
   \label{fig:RC}
\end{figure}
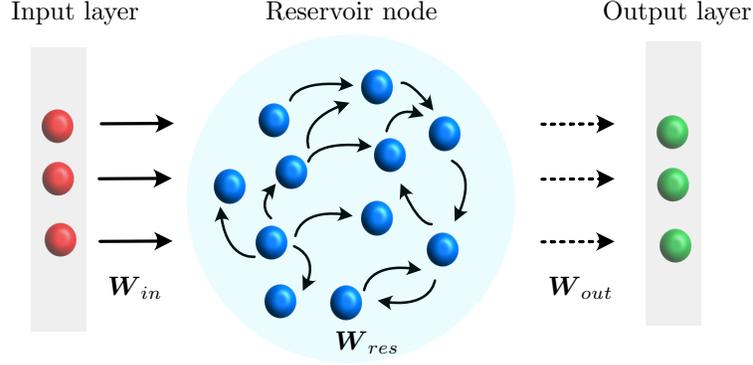

RC networks \cite {jaeger2002tutorial, lukovsevivcius2009reservoir, verstraeten2007experimental} are derived from recurrent neural network theory and 
map input signals into higher dimensional spaces 
through the dynamics of a non-linear system called 
reservoir. The idea in RC networks is to use recursive connections within neural networks to create complex dynamical systems. RC networks can be seen as a generalization of echo state networks (ESNs - see, e.g., \cite{jaeger2001echo})
and have the advantage of a straightforward and low-cost training process.
\fig{fig:RC} illustrates the main components of a typical ESN: i) an input layer, which in our case takes in the latent space representation of the FOM solution at time $t$ (i.e., $\bm{z}(t)$ in \eqref{equ_LS}); ii)
a hidden (between the visible input and output) %\textcolor{blue}{ In neural network community, the layer of neurons between input and outputs are typically called hidden layer, we can call it also Reservoir Layer }} 
layer of randomly connected neurons, i.e., the so-called reservoir node; and iii) an output layer, which returns an object belonging to a high dimensional space. The mapping to a richer space allows the RC framework to effectively capture the temporal dependency of the input signal. 

Let $N_h$ be the number of reservoir neurons. 
We denote with $\bm{h}(t) \in \mathbb{R}^{N_h}$  
% \anna{(per favore aggiungi le dimensioni)}
the reservoir state at time $t$ and with $\alpha \in [0, 1]$ the so-called leak rate. The evolution of the reservoir state is described by: 
\begin{equation}
\bm{h}(t^n) = (1 - \alpha) \bm{h}(t^{n-1}) + \alpha \tanh(\bm{W}_{\text{in}} \ \bm{z}(t^n) + \bm{W}_{\text{res}} \ \bm{h}(t^{n-1})),
\label{RC_state}
\end{equation}
where $\bm{W}_{\text{in}} \in \mathbb{R}^{N_h \times N_d}$ and $\bm{W}_{\text{res}} \in \mathbb{R}^{N_h \times N_h}$ 
% \anna{(anche qui)} 
are the input and reservoir weight matrices, respectively. See \fig{fig:RC}. Matrices $\bm{W}_{\text{in}}$ and $\bm{W}_{\text{res}}$ are initialized randomly and remain fixed throughout the entire process. 
Hyperparameter $\alpha$ controls the ``leaking'' of the previous state $\bm{h}(t^{n-1})$ into the current state $\bm{h}(t^n)$ and it determines how quickly the reservoir states 
are updated in response to new input $\bm{z}(t^n)$. One can say that $\alpha $ controls the balance between memory and reactivity in the reservoir. In fact, 
when $\alpha$ is close to one, the reservoir quickly adapts to the new input and the current reservoir state will reflect more the impact of the new input than the past information (long-term memory). %However,  high values of $\alpha$ increase the risk of forgetting past information.
On the other hand, for smaller values of $\alpha$, the reservoir state is updated more slowly and it retains 
more information from past states. This can help the network effectively memorize past information, which is beneficial for time series with long-term memory dependency.

The matrix of the output weights $\bm{W}_{\text{out}} \in \mathbb{R}^{N_d \times (N_d +N_h )} $ 
is learned during the training process using linear regression. This simplification aims to reduce the expensive computational 
cost of training to a linear regression task.
Despite it, RC remains effective in modeling and predicting complex dynamics, including chaotic time series \cite{lukovsevivcius2009reservoir}. To write down how 
$\bm{W}_\text{out}$ is obtained, let us introduce the following two matrices \cite{lukovsevivcius2009reservoir, trouvain2022create}:
\[
\mathbf{X} =
\begin{bmatrix}
    \bm{z}(t^0) & \bm{z}(t^1) & \dots & \bm{z}(t^{N_t - 1}) \\
    \bm{h}(t^0) & \bm{h}(t^1) & \dots & \bm{h}(t^{N_t - 1})
\end{bmatrix}, 
\quad
\mathbf{Y} =
\begin{bmatrix}
    \bm{z}(t^1) & \bm{z}(t^2) & \dots & \bm{z}(t^{N_t})
\end{bmatrix}.
\]
% \[
% \bm{X} =
% \begin{bmatrix}
%     \bm{z}(0) & \bm{h}(0) \\
%     \bm{z}(t^1) & \bm{h}(t^1) \\
%     \vdots & \vdots  \\
%     \bm{z}(t^{n-1}) & \bm{h}(t^{n-1})
% \end{bmatrix}, 
% \quad
% \bm{Y} =
% \begin{bmatrix}
%     \bm{z}(t^1) \\
%     \bm{z}(t^2) \\
%     \vdots \\
%     \bm{z}(t^n) 
% \end{bmatrix} .
% \]
%\textcolor{blue}{Ci sono due modi da scrivere X e Y, quale preferite voi?  io preferirei la versione piu compressa, version 2 }
% \textcolor{blue}{version 2}
% \[
% \mathbf{X} = 
% \begin{bmatrix}
%     \mathbf{z}(t^n) \\ \mathbf{h}(t^n)
% \end{bmatrix}_{n=0}^{N_t -1}, 
% \quad
% \mathbf{Y} = \left[ \mathbf{z}(t^{n}) \right]_{n=1}^{N_t}
% \]
%\anna{Arash, puoi sistemare l'ultima entrata di Y? Mi pare l'indice sia sbagliato.} \textcolor{blue}{ho corretto}
Matrix $\bm{X} \in \mathbb{R}^{(N_d +N_h ) \times N_t }$ is the so-called design matrix and matrix
$\bm{Y} \in \mathbb{R}^{N_d \times N_t }$ represents the target output. We recall that {$N_t$ is the number of training samples}.
%\anna{(diamo la taglia delle matrici per chiarire)}. \anna{Possiamo spiegare il ruolo della colonna di 1?} 
Then, we write:
\begin{align}
    \bm{W}_\text{out} = \bm{Y}\bm{X}^T(\bm{X} \bm{X}^T + \lambda \bm{I})^{-1},
\end{align}
where $\lambda$ is a regularization parameter that stabilizes the solution and reduces overfitting and $\bm{I} \in \mathbb{R}^{(N_d+N_h) \times (N_d+N_h)}$ 
% \anna{mettiamo le dimensioni anche qui per chiarezza}
is the identity matrix.

Given input $\bm{z}(t^n)$, the reservoir state $\bm{h}(t^n)$ is computed using \eqref{RC_state} and
then
the future latent space representation $\bm{z}(t^{n+1})$ is  predicted using the trained output weight matrix $\bm{W}_\text{out}$ as follows:
\begin{equation}
    \bm{z}(t^{n+1}) = \bm{W}_{\text{out}} 
    \begin{bmatrix}
    \bm{z}(t^n) \\
    \bm{h}(t^n)
    \end{bmatrix}
\label{RC_pred}
\end{equation}

{Once $\bm{z}(t^{n+1})$ is predicted with the equation above, it is fed to the decoder
of the CAE or E-CAE network \eqref{equ_dec} to obtain the corresponding physical solution.}

\begin{rem}
{While this work focuses on predicting one variable of interest ($\theta '$), the proposed framework is not limited to only one physical variable and it can be extended to multi-variable predictions at the same time. The input structure of the CAE would have to be modified from a single channel to a multi-channel input,  where each channel includes different physical variable at the same time instance. The CAE will learn to encode the dataset to a joint latent space that captures the dynamics of multiple variable.
Once the joint latent space is obtained, the subsequent procedure remain unchanged. The RC framework is trained to predict the future dynamics of this latent space. Finally, the CAE decodes the predicted latent space to physical space where each channel represents a different predicted physical variable. %This makes the framework flexible for predicting the multiple coupled variables simultaneously. 
}
\end{rem}

% \anna{Voglio che diamo la taglia delle matrici per assicurarci che 
% $\bm{z}(t+1)$ qui sopra e' delle taglia giusta. Al momento non mi e' chiaro} \textcolor{blue}{ho calcolato, la dimesione di $z(t^{n+1})$ viene $N_d$ perche $W_{out }$ ha la dimesnione $N_d \times (N_d +N_h )$ e vettore ha la dimesnione $  (N_d +N_h ) \times 1$ , cosi la dimesnion di $z(t^{n+1})$  viene $N_d$. }

%\anna{Arash, non capisco una cosa pero': fino a qui abbiamo spiegato come 
%predire il latent space futuro, giusto? Dobbiamo completare dicendo
%come passiamo dal latent space a $\theta'$, che e' quello che vogliamo predire
%di fatto.}

\section{Numerical results}\label{sec:res}
% \graphicspath{{./RTB/}}

Our proposed ROM approach has been validated using two well-known atmospheric flow benchmarks where a neutrally stratified atmosphere with uniform background potential temperature is perturbed by a bubble of either warm air (2D rising thermal bubble benchmark \cite{ahmad2007euler,ahmad2018high,feng2021hybrid,marras2015stabilized,hajisharifi2024comparison}) or cold air (2D density current benchmark \cite{ahmad2007euler,carpenter1990application,giraldo2008study,marras2013variational,marras2015stabilized,straka1993numerical,hajisharifi2024comparison}). There exist several variations of these two benchmarks, with different geometry and/ or initial condition. We adopted the setting from  \cite{ahmad2007euler}
for the rising thermal bubble and the setting from \cite{carpenter1990application, straka1993numerical} for the density current. 
Through these benchmarks, we aim at 
checking the accuracy of our ROM technique in reconstructing the time evolution and predicting the future dynamics of the potential temperature perturbation.
Results for the rising thermal bubble are presented in Sec.~\ref{sec:RTB}, while the results for the density current are discussed in Sec.~\ref{sec:DC}.
Finally, in Sec.~\ref{sec:RTB3D} we consider a three-dimensional variant of the rising thermal bubble benchmark to show that the proposed ROM works as well in 3D.

\subsection{Rising thermal bubble}
\label{sec:RTB}
\graphicspath{{Image/RTB/}}

% \anna{Dobbiamo stare attenti a non copiare troppo dai nostry paper, altrimenti quando sottomettiamo il paper potremmo avere problemi per autoplagio.
% Se hai gia' cambiato il testo vecchio, siamo a posto :)}
% \textcolor{blue}{si, ho gia fatto paraphrasing e il testo e' diverse dal  primo paper}

We perturb a neutrally stratified atmosphere with a uniform background potential temperature of $\theta_0$ = 300 K with a circular bubble of warmer air within computational domain $\Omega = [0,5000] \times [0, 10000]$ m$^2$ in the $xz$-plane. 
Over the time interval of $(0, 1020]$ s, the warm bubble rises due to buoyancy forces and transforms into a mushroom-like structure under the effect of shear stress.

The following initial temperature is prescribed:
\begin{align}
\theta^0 = 300 + 2\left[ 1 - \frac{r}{r_0} \right] ~ \textrm{if $r\leq r_0=2000~\mathrm{m}$}, \quad\theta^0 = 300
~ \textrm{otherwise},
\label{warmEqn1}
\end{align}
where $(x_c,z_c) = (5000,2000)~\mathrm{m}$ represents the center and $r = \sqrt[]{(x-x_{c})^{2} + (z-z_{c})^{2}}$ defines the radius of the bubble \cite{ahmad2007euler,ahmad2018high}. 
The initial density and initial specific enthalpy are given by: 
\begin{align}
\rho^0 &= \frac{p_g}{R \theta_0} \left(\frac{p}{p_g}\right)^{c_{v}/c_p} \quad \text{with} \quad p = p_g \left( 1 - \frac{g z}{c_p \theta^0} \right)^{c_p/R}, \label{eq:rho_wb} \\
h^{0} &= c_p \theta^0 \left( \frac{p}{p_g} \right)^{\frac{R}{c_{p}}},
\label{eq:e0}
\end{align}
with $p_g = 10^5$ Pa, $c_p = R + c_v$, $c_v = 715.5$ J/(Kg K), and $R = 287$ J/(Kg K). 
The initial velocity field is zero everywhere.
Impenetrable, free-slip boundary conditions are imposed on all the boundaries. 

We consider a uniform structured mesh with a grid size of $h = \Delta x = \Delta z = 62.5$ m. The time step is $\Delta t = 0.1$ s. We set $\mu_a = 15$ and $Pr = 1$ in \eqref{eq:mom_LES}-\eqref{eq:ent_LES} following \cite{ahmad2007euler,girfoglio2023validationAIP}.

%\anna{Possiamo spostare la frase seguente alla sez. in cui spieghiamo il metodo. Possiamo dire
%li' che usiamo $\mu_a$ costante e poi scrivere le
%frase seguente.} \textcolor{blue}{ho spostato  }

The original database is constructed by collecting the computed  potential temperature perturbation $\theta^\prime$ every
second. This database is then divided into a training dataset and a validation dataset.  As usual, the training dataset is used to train the CAE and RC networks, while the validation dataset is used to validate the ability of the CAE to reconstruct the physical space from the corresponding latent space and to asses the prediction accuracy of the RC network.

The CAE network we use for this benchmark has $L = 4$
convolutional layers with $N_f^1 = 256$, $N_f^2 = 128$, $N_f^3 = 64$, and $N_f^4 = 32$ and it gives
a latent space representation with $N_d = 4$.
See \fig{fig:CAE}. For the RC network, we choose
$N_h =400$ reservoir neurons and set the  leak rate to $\alpha =0.0095$ and 
the regularization parameter to $\lambda =0.004$ %\textcolor{blue}{fatto}.

%\anna{Arash scrive un paragrafo su come usare la stessa architettura se uno volesse predire altre variabili, e.g., $\rho$ o $h$.}

% This extension makes the framework highly flexible, allowing it to model and predict multiple coupled physical variables simultaneously without altering the core methodology.

%\anna{Curiosita': se volessimo predire un'altra variabile scalare (per es., ), useresti la stessa architettura per CAE e RC? O pensi che andrebbe tutto rifatto da zero (intendo settare numero di layers, filtri etc.)?} \textcolor{blue}{si, si pou usare lo stesso CAE e RC, ma dobbimao mettere il nuovo parametero sulla nuova canale di CAE, ne parliamo nel nostro meeting}

\begin{figure}[htb!]
    %\centering \includegraphics[width=160mm,scale=0.6]{Latent_space_RTB_All.pdf}
    \centering \includegraphics[width=\textwidth]{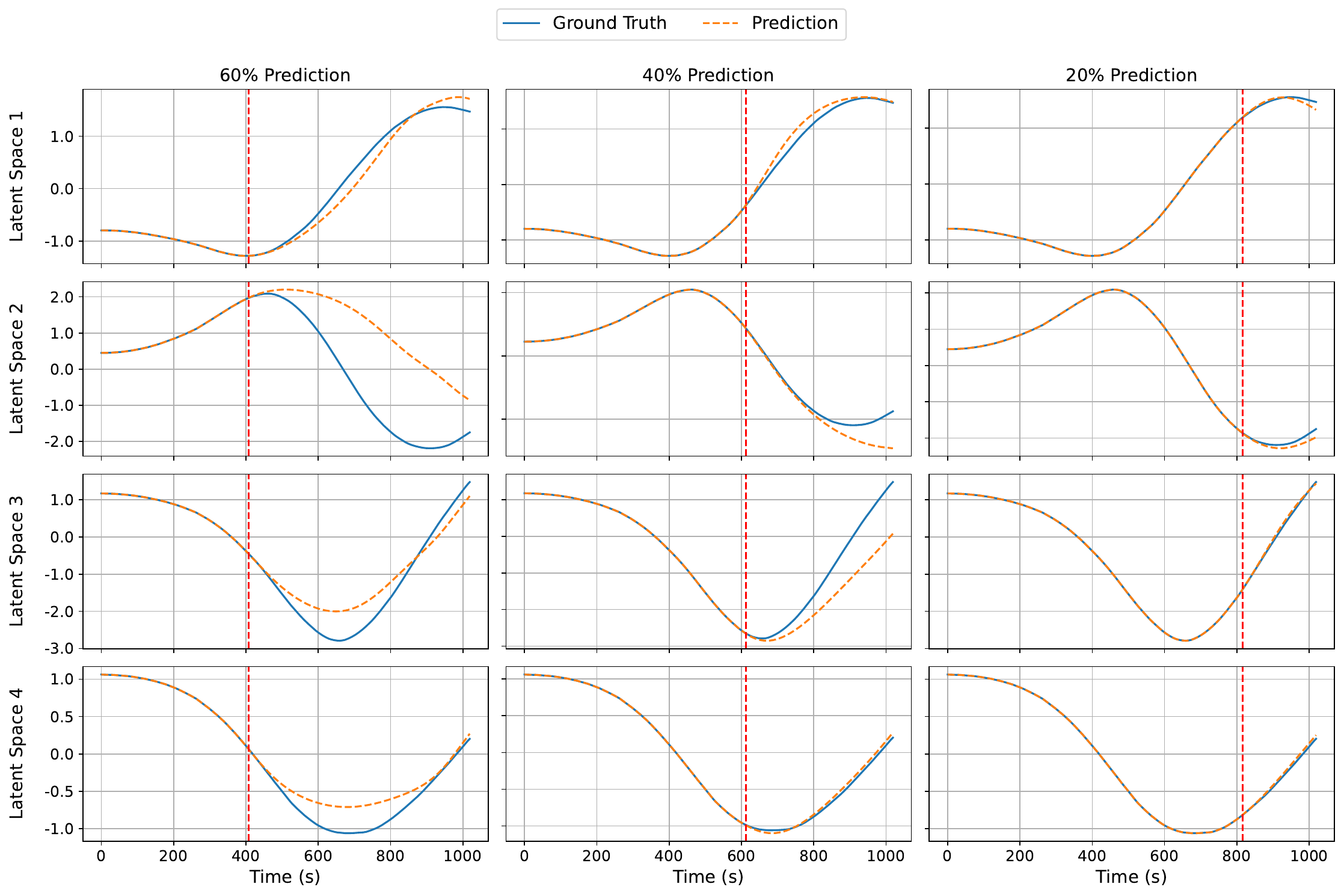}
        \caption{Rising thermal bubble: Time evolution of latent spaces for the different training-to-validation (T-to-V) dataset ratios listed in Tab. \ref{tab:sensitivity_analysis}.  The blue and orange dashed curves show the ground truth and the reconstructed/predicted latent spaces, respectively. The red dashed line indicates the time instance where prediction starts.}  
        \label{fig:LS_RTB}
\end{figure}

Let us analyze the model sensitivity to variations of the training-to-validation (T-to-V) dataset ratio. Tab.~\ref{tab:sensitivity_analysis} reports the three cases we considered with different T-to-V dataset ratios.
In all cases,  the snapshots corresponding to the first $N_t$ times are used to form the training set, with  $N_t = 408$ for the $40\% - 60\%$ ratio, $N_t = 612$ for the $60\% - 40\%$ ratio and $N_t = 816$ for the $80\%-20\%$ ratio. 
The remaining snapshots are reserved for validation. 
% The datasets are encoded  into 4 latent spaces using the CAE network depicted in \fig{fig:CAE} and 
\fig{fig:LS_RTB} shows the evolution of the CAE latent spaces for all cases reported in Tab.~\ref{tab:sensitivity_analysis}. The blue curve shows the latent space obtained from the CAE, denoted as ground truth, while the orange curve represent the latent spaces predicted by the RC network. We observe that the RC network does a poor job when only $40\%$ 
of the ground truth is used for training, due to insufficient data and a relatively long prediction window. The performance improves when we use $60\%$ of the database to predict the remaining $40 \%$, but the improvement becomes significant only when $80 \%$ of the database is used for training. In this last case, we obtain an excellent agreement between ground truth and reconstruction/prediction for all four latent spaces. 
%\anna{Altra curiosita': si puo' trovare una spiegazione per questi 4 valori? Se guardi bene, i primi 2 sono uno lo specchio dell'altro (piu' o meno) e gli ultimi 2 si assomigliano, con uno che scende al minimo piu' velocemente dell'altro.
%Qui l'immagine e' abbastanza semplice (una palla che sale e si deforma), per questo me lo chiedevo.} \textcolor{blue}{come viene l'andamneto del Latent space e' balck box per noi. Intendo che non abbimo nessun controllo sul andamaneto del latent space e viene queelo che viene. per questo non possiamo interpretare latent space. forse si non sosno 100 per cento sicuro, ma io non ho visto }

\begin{table}[htb!]
\centering
\begin{tabular}{|c|c|c|c|}
\hline
\rowcolor{gray!30} & \textbf{T-to-V dataset ratio } & \textbf{T Snapshots \#}& \textbf{V Snapshots \#} \\ \hline
\textbf{Case 1}      & 40\% - 60\%           & 408          & 612             \\ \hline
\textbf{Case 2}      & 60\% - 40\%             & 612          & 408             \\ \hline
\textbf{Case 3}      & 80\% - 20\%             & 816          & 204             \\ \hline
\end{tabular}
\caption{Rising thermal bubble: training-to-validation (T-to-V) dataset ratio, number of training snapshots and validation snapshots 
for the three cases used for the sensitivity analysis.}
\label{tab:sensitivity_analysis}
\end{table}

%This demonstrates the capability of our model to accurately predict future dynamics in the presence of more training data. 

%\anna{Arash aggiunge un commento sul fatto che il latent space cambia ogni volta che si fa il training e che deve essere il piu' smooth possible affinche' la nostra procedura funzioni bene.} 

We note that CAE networks introduce multiple sources of randomness, which are hard to control.
This makes it difficult to obtain an identical latent space representations across different training instances. Consequently, different latent space representations may be obtained each time the CAE is trained. It is important for this representation to be as smooth (see the example in \fig{fig:LS_RTB}) as possible for our ROM to work seamlessly, as irregularities could worsen the accuracy of our predictive model.

\begin{figure}[htb!]
    \centering \includegraphics[width=120mm,scale=1.0]{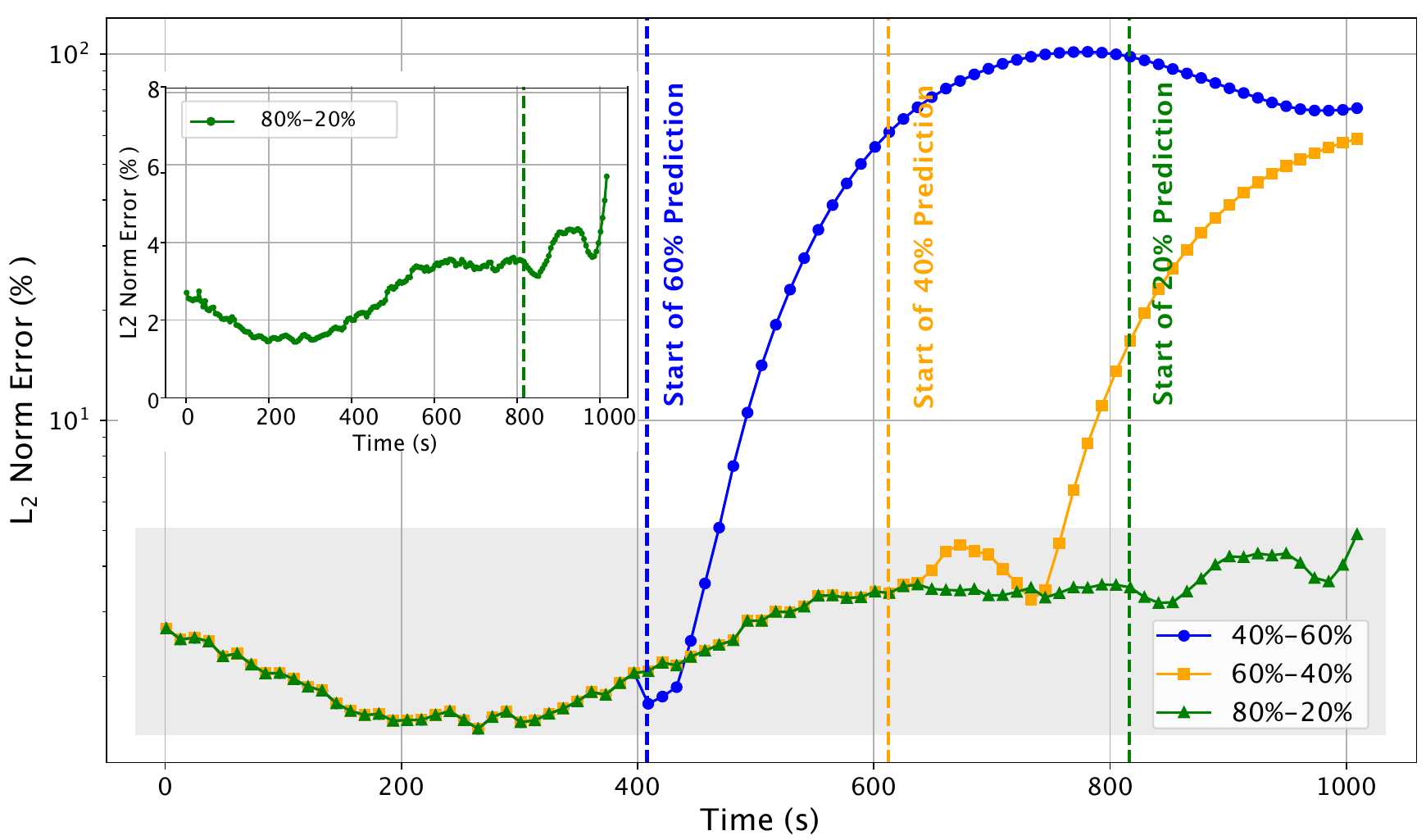}
        \caption{Rising thermal bubble: Time evolution  of error \eqref{eq:l2Error} for different training-to-validation (T-to-V) dataset ratios listed in Tab. \ref{tab:sensitivity_analysis}: $60\% - 40\%$ (blue curve), $40\% - 60\%$ (yellow curve) and $80\% - 20\%$ (green curve). The dashed vertical lines indicate the start of prediction phase for each case. The inset provides a zoomed-in view of green curve in normal scale.}  
        \label{fig:Rec_error_RTB_All}
\end{figure}

For a quantitative assessment of the reconstructions and predictions,
we compute the $L^2$ error between FOM and ROM solutions:
\begin{equation}
E_{\theta^\prime}(t) = 100 \cdot \dfrac{||\theta^\prime_h(t) - \theta^\prime_r(t)||_{L^2(\Omega)}}{||{\theta^\prime_h}(t)||_{L^2(\Omega)}}.
\label{eq:l2Error}
\end{equation}
%\anna{dobbiamo assicurarci di usare il pedice $h$
%quando parliamo della soluzione FOM nella discussione del metodo.} \textcolor{blue}{ok}
\fig{fig:Rec_error_RTB_All} shows error
\eqref{eq:l2Error} for the three T-to-V dataset ratios under consideration. In all cases, the reconstruction error (i.e., the error within the training time interval) remains low, demonstrating that the CAE is capable of accurately reconstructing the FOM solution from the latent spaces. Additionally, this suggests that the CAE can perform its reconstruction task effectively even when trained on only $40 \%$ of the database. 
However, the prediction error varies considerably in the three cases. 
For the case 1 (blue curve), the error increases rapidly after the beginning of prediction phase. This behaviour points to a clear insufficiense of training data.  For case 2 (orange curve), the error remains lower than $10\%$ for the first $200$ s of prediction, i.e., until $t =800$ s.
However, after $t =800$ s the prediction worsens quickly. In case 3 (green curve), we obtain a low error ($< 6\%$) throughout the entire prediction phase, which lasts about 200 s.
However, even in this third case, we see a sharp
increase in the error towards the end of the prediction phase, indicating that one could probably not expect the error to remain low for much longer. 
\fig{fig:Rec_error_RTB_All} suggests that the RC network can accurately predict the future dynamic when trained on a sufficiently large amount of training data. 
Hence, from now on we stick to the 
$80\%$-$20\%$ T-to-V dataset ratio to 
analyze the performance of our ROM for this benchmark.

% \begin{figure}[t]
%     \centering
%     \begin{tikzpicture}
%         % Main figure
%         \node[anchor=south west,inner sep=0] (main) at (0,0) 
%             {\includegraphics[width=1.0\textwidth]{Rec_Err_RTB_All.pdf}};

%         \put(14,40.5){\tiny{\textcolor{white}{$t$=400 s}}}
%         % Zoomed View
%         \node[anchor=south west,inner sep=0] (inset) at (1.3,3.5)  
%             {\includegraphics[width=0.35\textwidth]{Rec_Error_RTB.pdf}};
        
%         % Draw a rectangle on the main figure
%         \draw[red,thick] (1.75,1.3) rectangle (14.5,3.5); 
    
%         % Draw an arrow 
%         \draw[->,thick,red] (6,2) -- (4,4.5); 
%     \end{tikzpicture}

%     \caption{Rising thermal bubble: Time evolution  of the $L^2$ error \eqref{eq:l2Error} between FOM and ROM solutions for different training-to-validation (T-to-V) dataset ratios listed in Tab. \ref{tab:sensitivity_analysis}: $60\% - 40\%$ (blue curves), $40\% - 60\%$ (yellow curves) and $60\% - 40\%$ (green curve). The dashed vertical lines indicates the start of prediction phase for each case. The inset provides a zoomed-in view of green curve in normal scale.}
%     \label{fig:Rec_error_RTB_All}
% \end{figure}

\begin{figure}[htb!]
    \vspace{1cm}
    %\begin{subfigure}{\linewidth} 
    \centering
    \begin{overpic}[width=0.56\textwidth, grid=false]{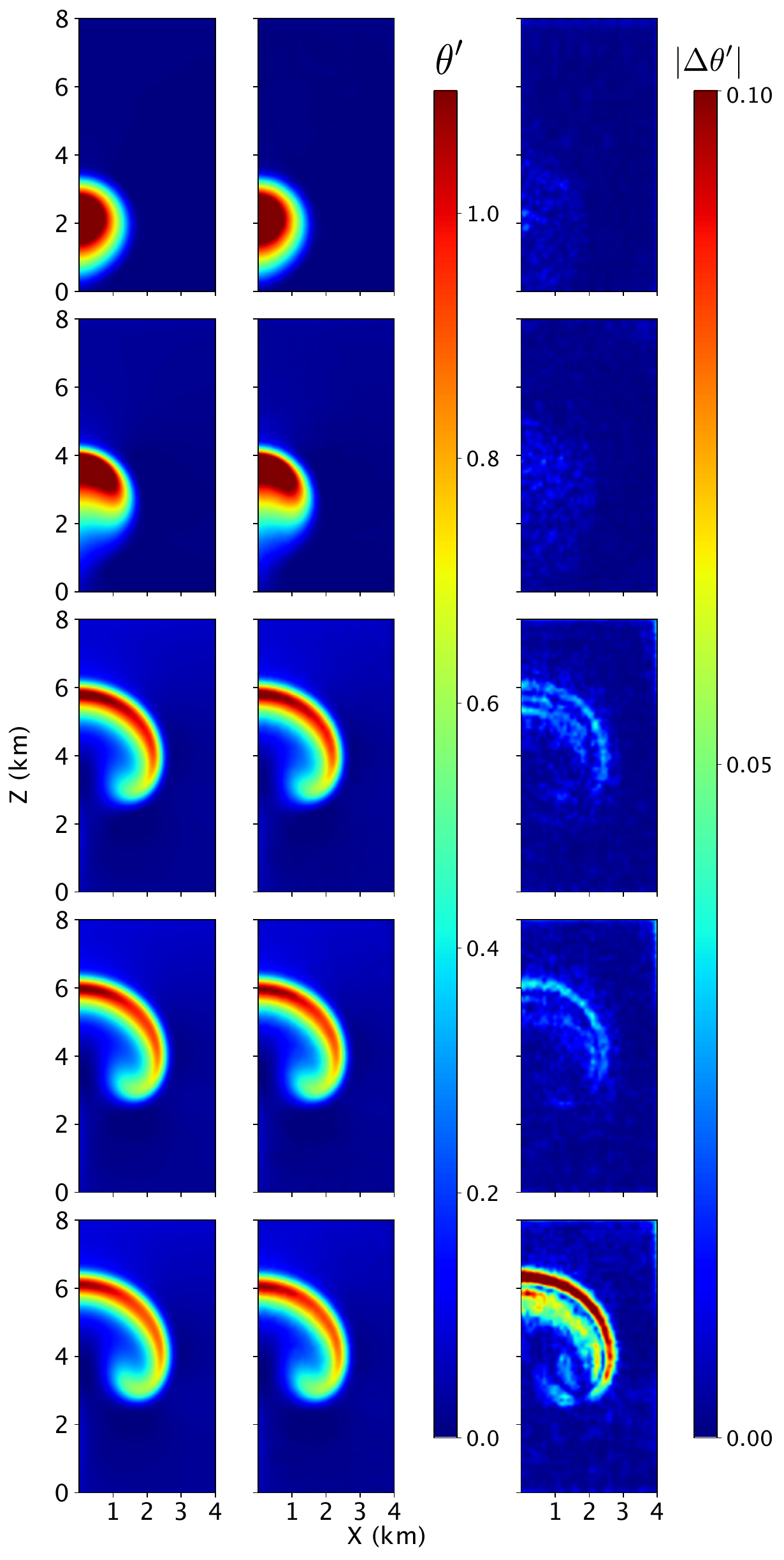}  
         \put(6.8,100.5){ROM}
         \put(18.2,100.5){FOM}
         \put(36.2,100.5){$|\Delta \theta ' |$}
         \put(6.4,97.5){\tiny{\textcolor{white}{$t$=255 s}}}
         \put(6.4,78.3){\tiny{\textcolor{white}{$t$=505 s}}}
         \put(6.4,59.0){\tiny{\textcolor{white}{$t$=930 s}}}
         \put(6.4,39.6){\tiny{\textcolor{white}{$t$=980 s}}}
         \put(6.4,20.5){\tiny{\textcolor{white}{$t$=1020 s}}}
         
          \put(19,97.5){\tiny{\textcolor{white}{$t$=255 s}}}
         \put(19,78.3){\tiny{\textcolor{white}{$t$=505 s}}}
         \put(19,59.0){\tiny{\textcolor{white}{$t$=930 s}}}
         \put(19,39.6){\tiny{\textcolor{white}{$t$=980 s}}}
         \put(18,20.5){\tiny{\textcolor{white}{$t$=1020 s}}}

          \put(35.8,97.5){\tiny{\textcolor{white}{$t$=255 s}}}
         \put(35.8,78.3){\tiny{\textcolor{white}{$t$=505 s}}}
         \put(35.8,59.0){\tiny{\textcolor{white}{$t$=930 s}}}
         \put(35.8,39.6){\tiny{\textcolor{white}{$t$=980 s}}}
         \put(34.8,20.5){\tiny{\textcolor{white}{$t$=1020 s}}}

      \end{overpic}
      \captionsetup{list=no} 
  \caption{Rising thermal bubble: Comparison of the evolution of $\theta'$
   given by the ROM (first column) and the FOM (second column) and their differences in absolute value (third column). %\anna{Arash, puoi mettere il valore assoluto per $\Delta \theta'$ cosi' da essere piu' precisi?} 
  %\textcolor{blue}{Fatto} 
  }
    % \anna{E' una mia impressione o le figure sono distorte?} \textcolor{blue}{E' risoloto } \anna{Io le vedo ancora distorte. Confronta con la Fig. 3 del nostro paper precedente per vedere la differenza. Come mai la differenza ha la colorbar che va fino a 0.6? NOn riesci a vedere rosso in nessuna delle figure, forse va tarata meglio?}}  
   \label{fig:Vis_RTB}
\end{figure}

\fig{fig:Vis_RTB} displays a qualitative comparison between ROM and FOM solutions at five times. Out of these five time instances, the first two (i.e., $t = 255$ s and $t = 505$ s) are part of the training dataset and are meant to demonstrate the ability of the CAE to reconstruct the solution from the latent spaces. For these 2 times, we see excellent agreement between ROM and FOM solutions, as expected from Fig.~\ref{fig:Rec_error_RTB_All}.
The remaining three time instances (i.e., $t = 930,980, 1020$ s) belong to the validation dataset and are used to assess the accuracy of the ROM to predict the future state of the system.
We observe that the ROM is able to predict the FOM results rather accurately, although at the final time $t = 1020$ s the agreement is less satisfactory.
This was expected since, as shown in the sub-panel of \fig{fig:Rec_error_RTB_All}, 
the $L^2$ error between FOM and ROM solutions for case 3 ($80\%$-$20\%$ T-to-V ratio) reaches 
$6 \%$ at the end of simulation.   
The third column in \fig{fig:Vis_RTB} shows the 
difference between ROM and FOM solutions in 
absolute value with a color bar narrowed to
interval $[0, 0.1]$, which corresponds to $10\%$ of the maximum value, to better visualize the differences.
%\anna{Comment on the difference once the colobar is adjusted}

% \anna{E' meglio se spieghiamo di piu' 
% quest'ultimo punto. Non sono sicura di 
% aver capito cosa si intende.}\textcolor{blue}{agiungo anche qui la differenza tra ROM e FOM per essere piu chiaro, ora e' piu chiaro  }

% \anna{Questo grafico e' lo stesso della curva verde in Fig.~2,no? In caso, non lo ripeterei.} \textcolor{blue}{ho tolto la figura}

% \begin{figure}[htb!]
%     \centering \includegraphics[width=100mm,scale=0.6]{Rec_Error_RTB.pdf}
%         \caption{Rising thermal bubble:  Reconstruction Error. }   
%         \label{fig:Rec_error_RTB}
% \end{figure}

All the simulations were run on a 11th Gen Intel(R) Core(TM) i7-11700 @ 2.50GHz system with 32 GB RAM. On this machine, a FOM simulation takes  $65$ s. {The training for the CAE and RC framework, which was conducted on Google Colab with A100 GPUs, takes approximately $1500$ s.  The time required by the ROM to predict the latent space and reconstruct the physical solution after training phase is around less than 1 s.}

% the time need to predict the future latent spaces and map them the into physical space is around $0.8$ s.  

%\textcolor{blue}{Abbiamo bisgono di aggiungere questa parte per ogni benchmar ? perche la simulzione 3D era piu lungo del 2D}

\subsection{Density Current}
\label{sec:DC}
\graphicspath{{Image/DC/}}

We perturb a neutrally stratified atmosphere with a uniform background potential temperature of $\theta_0$ = 300 K with a circular bubble of cold air within computational domain  $\Omega = [0, 25600] \times [0, 6400]$ m$^2$ in the $xz$-plane. During the time 
interval of interest $(0, 900]$ s, the cold air experiences two primary motions: initially, the negative buoyancy generates vertical motion driving the cold air toward the ground. Upon reaching the ground, the air forms a cold front that propagates through a dominant horizontal motion, during which a complex structure with multiple vortices is created.

The initial temperature field is given by: %\anna{usiamo lo stesso stile di (19)} \textcolor{blue}{fatto}
\begin{align}
    \theta_0 = 300 - \theta_s [1+ \cos(\pi r)], \quad \mathrm{if} \ r\le1, \quad \mathrm{otherwise} \ \theta_0 = 300, 
    \label{eq:initila_pot_temp_DC}
\end{align}
% \begin{equation}
% \theta^0 = 300 + 2\left[ 1 - \frac{r}{r_0} \right] ~ \textrm{if $r\leq r_0=2000~\mathrm{m}$}, \quad\theta^0 = 300
% ~ \textrm{otherwise},
% \label{warmEqn1}
% \end{equation}
where \(\theta_s\) represents the semi-amplitude of the initial temperature perturbation and it is set to 7.5 \cite{ahmad2007euler,carpenter1990application,giraldo2008study,marras2013variational,marras2015stabilized,straka1993numerical}. The parameter \(r\) represents a normalized radial distance, calculated as:
\begin{equation}
    r = \sqrt{\left(\frac{x - x_c}{x_r}\right)^2 + \left(\frac{z - z_c}{z_r}\right)^2}.
\end{equation}
with $(x_r,z_r) = (4000, 2000) ~\mathrm{m}$ and $(x_c, z_c) = (0, 3000) ~\mathrm{m}$.   
The initial density and specific enthalpy are given by \eqref{eq:rho_wb} and \eqref{eq:e0}, respectively, while the velocity field is initialized to zero across the domain. Impenetrable, free-slip boundary conditions are imposed on all walls. 

We consider a structured grid with a uniform spacing $h = \Delta x = \Delta z = 100$ m and time step is $\Delta t = 0.1$ s. We set $\mu_a = 75$ and $Pr = 1$ in \eqref{eq:mom_LES}-\eqref{eq:ent_LES} to stabilize the numerical solution \cite{ahmad2007euler,straka1993numerical}. To evaluate the accuracy of our ROM model, 
we sample the $\theta^\prime$ field at the frequency of 1 second, resulting in a dataset of 900 snapshots. Based on the sensitivity analysis of the training-to-validation dataset ratio conducted for the previous benchmark, we use a ratio of  80$\%$-20$\%$ (720 - 180 snapshots) to generate the training and validation datasets. 

The snapshots are encoded into $N_d = 4$ latent spaces through a CAE/E-CAE and a
RC network is subsequently trained on the obtained latent spaces to predict future. The CAE and E-CAE networks for this benchmark consist of $L = 4$ convolutional layers with $N_f^1 = 512$, $N_f^2 = 256$, $N_f^3 = 128$, and $N_f^4 = 64$. The number of filters per layer is increased with respect to the benchmark in the previous section because the density current benchmark features a more complex dynamics. 
%\anna{Arash, il motivo per cui questi valori di $N_f$ sono piu' elevati che per la RTB e' perche' le immagini sono piu' complesse?} \textcolor{blue}{Esatto, la dinamica qui e' piu complessa e abbimao bisogno di piu filter per catturare la dinamica della problema}
For the E-CAE network, we set $n_f=3$.
See \fig{fig:E-CAE}.
For the RC network, we choose
$N_h = 1000$ reservoir neurons and set the  leak rate to $\alpha = 0.0022 $ and 
the regularization parameter to $\lambda = 0.0022$. %\textcolor{blue}{fatto}

\begin{figure}[htb!]
    %\centering \includegraphics[width=100mm,scale=0.6]{Rec_Error_Comparison_CAEs.pdf}
    \centering \includegraphics[width=0.5\textwidth]{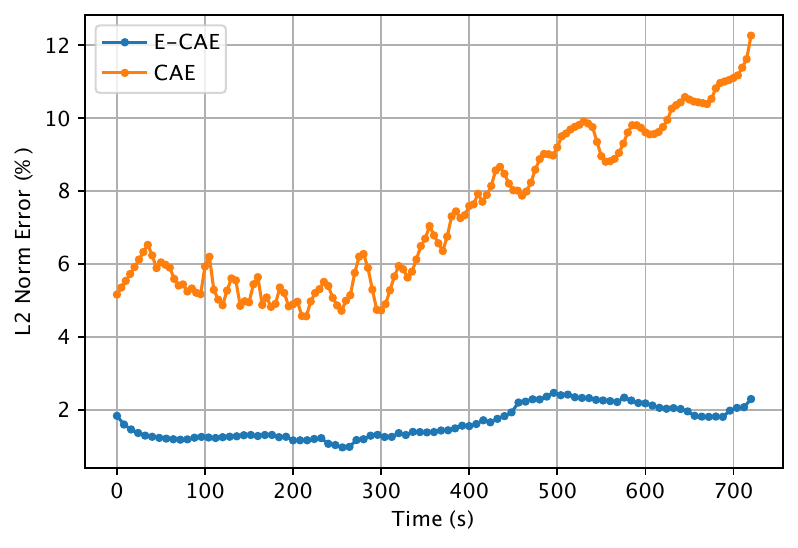}
        \caption{Density current:  Time evolution  of the $L^2$ error \eqref{eq:l2Error} between FOM and ROM solutions for CAE (orange curve) and E-CAE (blue curve). }   
        \label{fig:Rec_error_DC_2CAEs}
\end{figure}

\fig{fig:Rec_error_DC_2CAEs} shows 
the evolution of error \eqref{eq:l2Error} in the reconstruction phase $(0,720]$ s
for two convolutional autoencoder networks: CAE (orange curve) and E-CAE (blue curve). 
While the CAE network depicted in \fig{fig:CAE}
worked well for the simpler dynamic of 
the warm bubble benchmark, it struggles to capture the more complex dynamic of this benchmark. The error of the CAE increases over time and reaches around 12$\%$ by time
$t = 720$ s. This limitation is due to the insufficient depth of the network which constrains its ability to effectively learn the highly non-linear dynamics arising in this benchmark. 
From \fig{fig:Rec_error_DC_2CAEs}, it is evident that the department-wise strategy of E-CAE 
illustrated in \fig{fig:E-CAE}
%\anna{(il 2 e' un typo?)} \textcolor{blue}{Si, e' un typo, risolto} 
significantly outperforms a standard CAE: error \eqref{eq:l2Error} remains below $2\%$ for most of the time interval $(0,720]$ s, achieving $75 \%$ lower reconstruction error and fewer oscillations than the CAE. 
This suggests that the additional layers can successfully capture the complex dynamics, enabling the network to accurately reconstruct the physical field from its latent spaces. Hence, from the rest of the section
we will present results obtained with the E-CAE.

\begin{figure}[t]
    \centering 
    \includegraphics[width=0.7\textwidth]{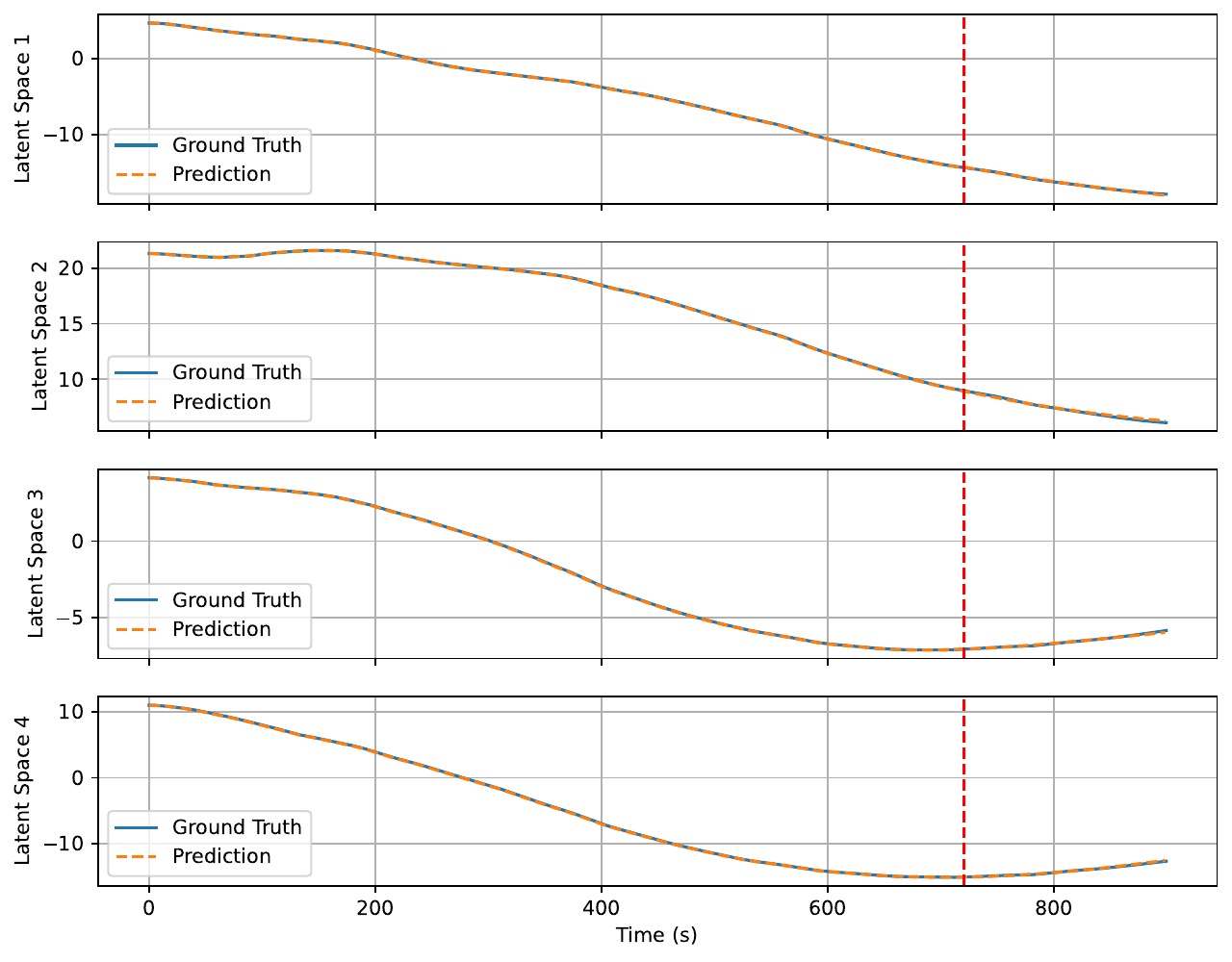}
        \caption{Density current: Time evolution of latent spaces obtained using E-CAE (blue curve) and reconstructed/predicted by RC framework (orange curve).  The red dashed line shows the beginning of prediction phase.}  
        \label{fig:LS_DC}
\end{figure}

\fig{fig:LS_DC} illustrates the latent spaces obtained through E-CAE (ground truth) and the latent spaces {reconstructed/predicted using RC network (prediction)}. %\anna{Mi sembra che usare prediction sia un po' fuorviante perche' e' reconstruction and prediction, no? Infatti anche la figura e' un po' ambigua perche' nella legenda c'e' "prediction" e poi diciamo che la linea rossa indica l'inizio della prediction phase.}. \textcolor{blue}{si hai ragione, ho modificato il testo, come ti sembara ora?}
The excellent agreement we observe in all four panels in \fig{fig:LS_DC} shows the capability of the RC network to effectively learn the dynamics of compressed data representation and accurately predict their future evolution.
%\anna{Anche qui mi chiedo il significato delle 4 curve. Le ultime 2 sono identiche, mi pare, solo shiftate.} %\textcolor{blue}{No, non sono uguali, sono diversi. Se sei d'accordo, Ne parliamo nella nostra riunione}

\begin{figure}[htb!]
    \centering \includegraphics[width=.5\textwidth]{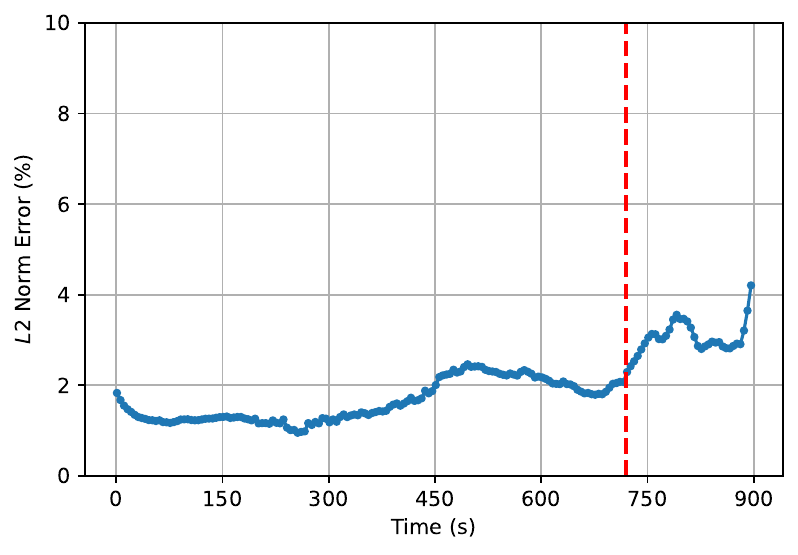}
    %\centering \includegraphics[width=100mm,scale=0.3]{Rec_Error_DC.pdf}
        \caption{Density current: Time evolution  of the $L^2$ error \eqref{eq:l2Error} between FOM and ROM solutions. The red dashed line marks the beginning of prediction. }   
        \label{fig:Rec_error_DC}
\end{figure}

\fig{fig:Rec_error_DC} shows the reconstruction error of the E-CAE till
$t= 700$ s and the combined error (the sum of the reconstruction and prediction errors) from $t=700$ s to $t=900$ s. The reconstruction error remains stable at around 2$\%$ throughout the training phase while the combined error gradually increases to 4$\%$, which is still a  rather remarkable performance of the RC network in predicting the future state of latent spaces.
However, once again we observe a sharp rise
in the error towards the end of the prediction 
phase, which indicates the accuracy in the prediction might not decline quickly past 
$t = 900$ s.

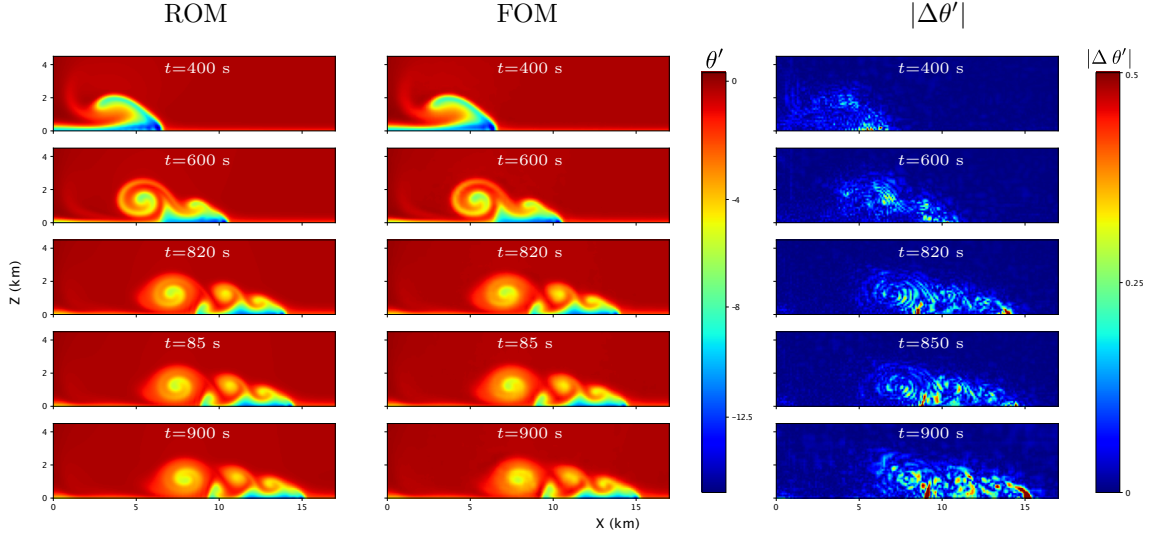
\begin{figure}[htb!]
    \vspace{1cm}
    %\begin{subfigure}{\linewidth} 
    \centering
    \begin{overpic}[width=1\textwidth, grid=false]{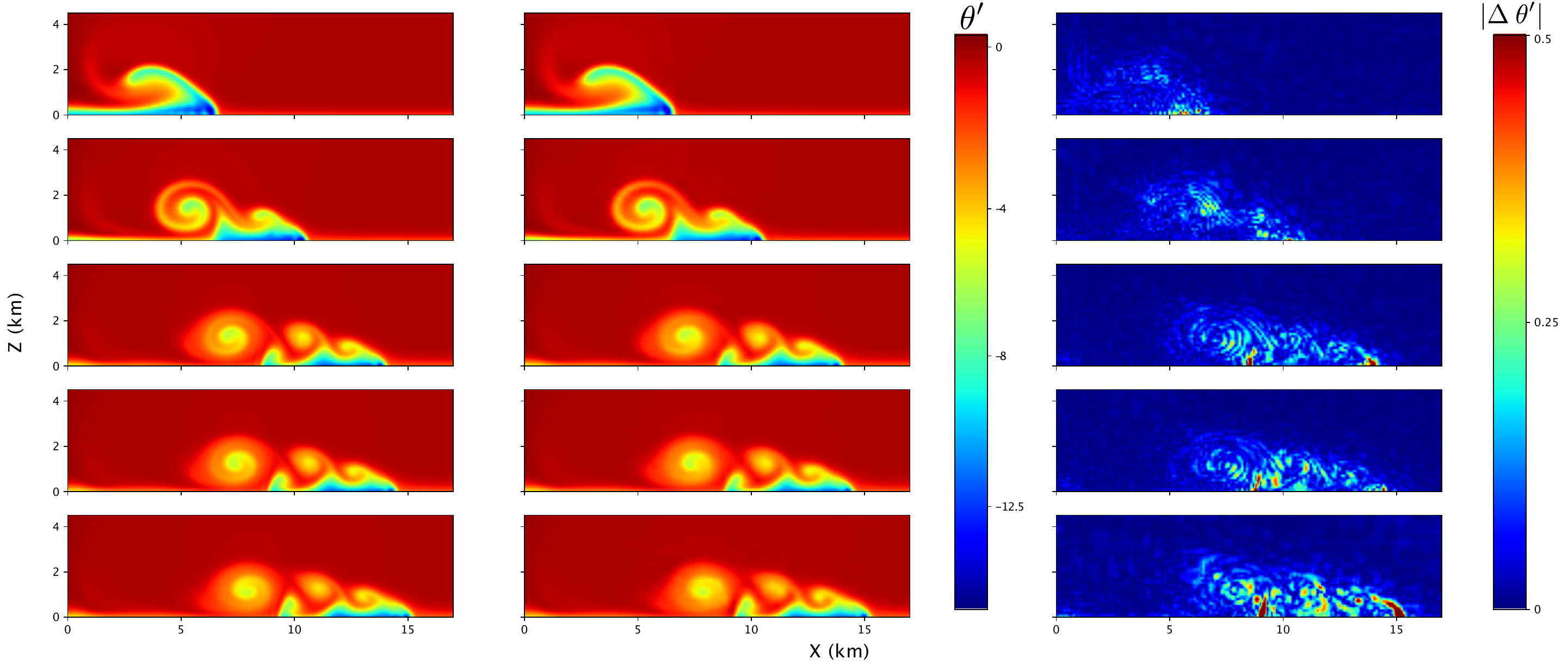}  
          \put(14,45){ROM}
            \put(43,45){FOM}
            \put(79,45){$| \Delta {\theta} ^\prime |$ }
            \put(14,40.5){\tiny{\textcolor{white}{$t$=400 s}}}
            \put(14,32.5){\tiny{\textcolor{white}{$t$=600 s}}}
            \put(14,24.5){\tiny{\textcolor{white}{$t$=820 s}}}
            \put(14,16.5){\tiny{\textcolor{white}{$t$=85 s}}}
            \put(14,8.7){\tiny{\textcolor{white}{$t$=900 s}}}

            \put(43,40.5){\tiny{\textcolor{white}{$t$=400 s}}}
            \put(43,32.5){\tiny{\textcolor{white}{$t$=600 s}}}
            \put(43,24.5){\tiny{\textcolor{white}{$t$=820 s}}}
            \put(43,16.5){\tiny{\textcolor{white}{$t$=85 s}}}
            \put(43,8.7){\tiny{\textcolor{white}{$t$=900 s}}}

            \put(78,40.5){\tiny{\textcolor{white}{$t$=400 s}}}
            \put(78,32.5){\tiny{\textcolor{white}{$t$=600 s}}}
            \put(78,24.5){\tiny{\textcolor{white}{$t$=820 s}}}
            \put(78,16.5){\tiny{\textcolor{white}{$t$=850 s}}}
            \put(78,8.7){\tiny{\textcolor{white}{$t$=900 s}}}

      \end{overpic}
      \captionsetup{list=no} 
  \caption{Density current: comparison of the evolution of $\theta'$ given by the ROM (first column) and the FOM (second column) and their differences in absolute value (third column). %\anna{Un paio di commenti: 1) non superare 1textwidth perche' data' problemi in fase di sottomission, 2) cambiamo il massimo di $\Delta \theta'$ anche qui per vedere delle regioni in rosso come abbiamo fatto in fig. 6 e 3)  aggiungiamo il valore assoluto per $\Delta \theta'$.}\textcolor{blue}{Fatto}
  }
   \label{fig:Vis_DC}
\end{figure}

\fig{fig:Vis_DC} displays a qualitative comparison between ROM and FOM solutions, 
with their difference in absolute value, for 5 time instances. 
Like in the case of \fig{fig:Vis_RTB}, 
the first two time instances (i.e., $t = 400, 600$ s) 
belong to the training phase and demonstrate the accuracy of the E-CAE. 
The close agreement that we see in the first two rows of \fig{fig:Vis_DC}
reflects the small $L^2$ error seen in \fig{fig:Rec_error_DC} during
the training phase.  
The last three rows in \fig{fig:Vis_DC}, corresponding to
$t = 820,850,900$ s, show the accuracy of E-CAE and RC networks in predicting 
the future state of the system. The small mismatch in predicting the propagation of 
the cold front and the multi-rotor flow structure demonstrate the 
accuracy of our approach. 
%\anna{Arash, puoi aggiungere anche qui un commento
%sulla color bar, che abbiamo forzato.}
To improve the visualization of the difference between the ROM and FOM solutions in absolute value, we
restricted the  color bar for the third column of \fig{fig:Vis_DC} to $[0, 0.5]$, corresponding to $4\%$ of the maximum value. 

{The FOM simulation on the same local hardware mentioned at the end of Sec.~\ref{sec:RTB} takes $286$ s. The training of the CAE and RC framework takes approximately $4100$ s. Then, the time needed for the ROM prediction is around $2$ s.}

%\anna{Arash aggiunge un commento sui tempi computazionali.}

\subsection{3D Rising Thermal Bubble}
\label{sec:RTB3D}
\graphicspath{{Image/RTB3D/}}

To demonstrate that the applicability of the ROM under consideration is not limited to 2D problems, we consider a tree-dimensional variant \cite{marras2013variational,hajisharifi2024comparison} of the rising thermal bubble benchmark  in Sec.~\ref{sec:RTB}.

\begin{figure}[]
    \vspace{1cm}
    %\begin{subfigure}{\linewidth} 
    \centering
    \begin{overpic}[width=0.72\textwidth, grid=false]{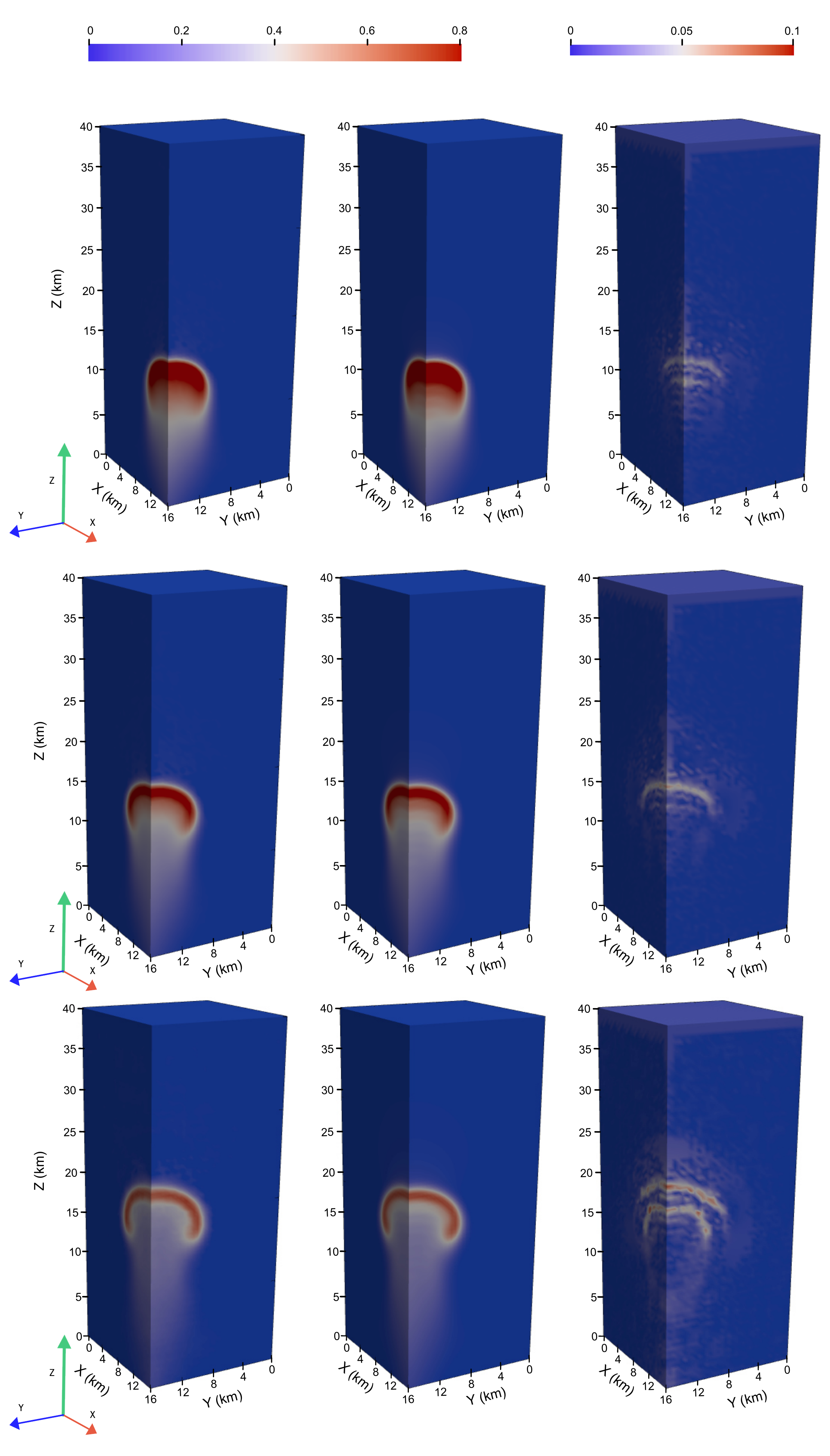}  
          
          \put(45,99){$| \Delta {\theta} ^\prime |$ }
          \put(18.5,99){$ {\theta} ^\prime $ }
    
          \put(11,93){ROM}
          \put(29,93){FOM}
          \put(47,93){$| \Delta {\theta} ^\prime |$ }
            % \put(79,45){$\Delta {\theta} ^\prime$ }
            \put(12.5,88){\tiny{\textcolor{white}{$t$=340 s}}}
            \put(30.,88){\tiny{\textcolor{white}{$t$=340 s}}}
            \put(48,88){\tiny{\textcolor{white}{$t$=340 s}}}

            \put(12.5,57){\tiny{\textcolor{white}{$t$=420 s}}}
            \put(30,57){\tiny{\textcolor{white}{$t$=420 s}}}
            \put(48,57){\tiny{\textcolor{white}{$t$=420 s}}}

            \put(12.5,27){\tiny{\textcolor{white}{$t$=500 s}}}
            \put(30,27){\tiny{\textcolor{white}{$t$=500 s}}}
            \put(48,27){\tiny{\textcolor{white}{$t$=500 s}}}

      \end{overpic}
      \captionsetup{list=no} 
  \caption{3D rising thermal bubble: evolution of $\theta'$ given by the ROM (left column) and the FOM (right column), and their differences in absolute value (third column).}
   \label{fig:Vis_RTB3D}
\end{figure}
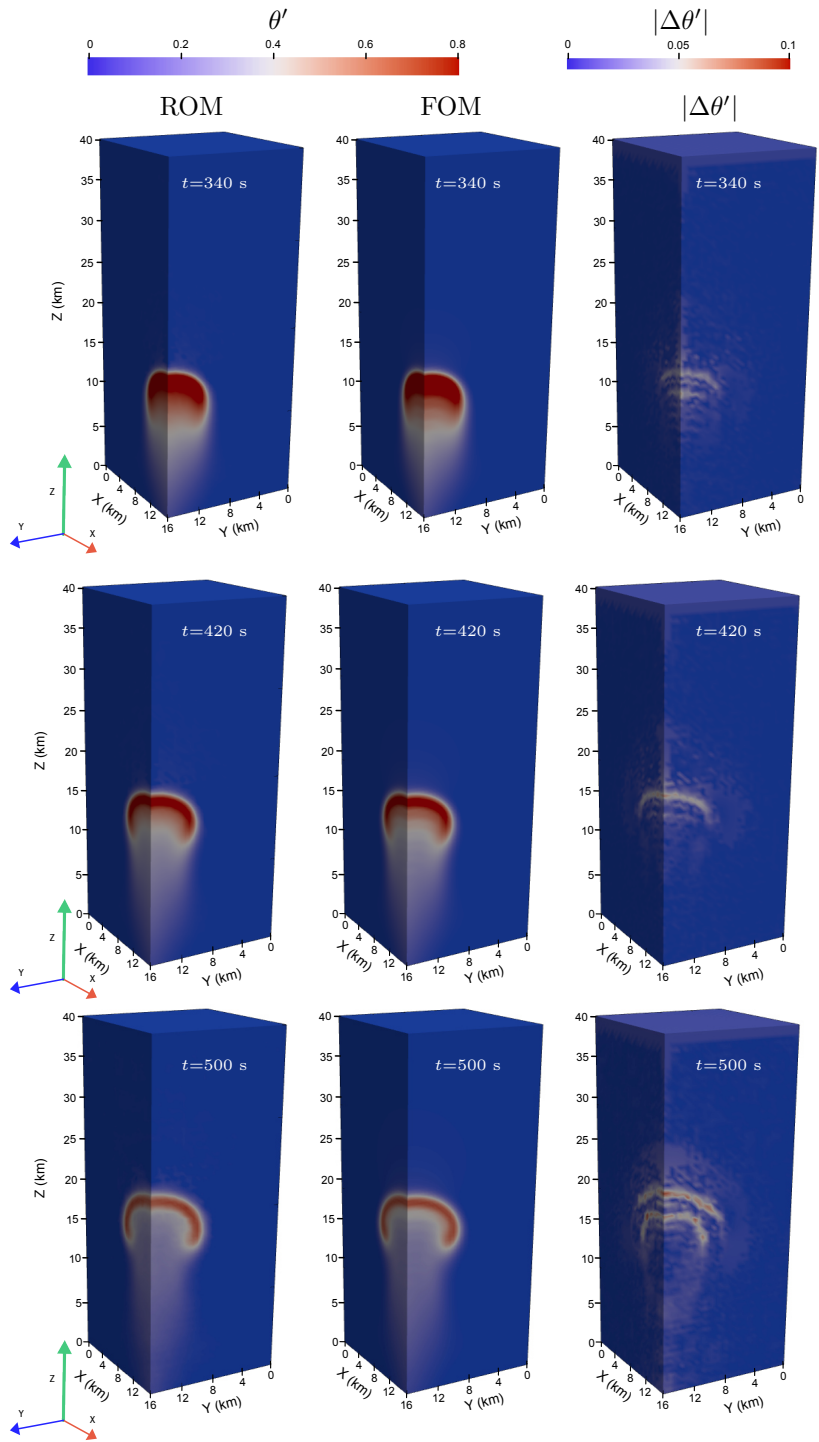

The computational domain is a parallelepiped with dimensions $\Omega = [0,1600] \times [0,1600] \times [0,4000] \, \mathrm{m}^3$. Similar to the 2D version of the benchmark, the system starts from a neutrally stratified atmosphere with uniform background potential temperature at
300 K, with a disturbance in the form of
a spherical bubble of warmer air. This initial temperature field is given by:
\begin{equation}
\theta^0 = 300 + 2\left[ 1 - \frac{r}{r_0} \right] ~ \textrm{if $r \leq r_0 = 500 \, \mathrm{m}$}, \quad \theta^0 = 300 ~ \textrm{otherwise},
\end{equation}
with $r = \sqrt[]{(x - x_c)^2 + (y - y_c)^2 + (z - z_c)^2}$
and $(x_c, y_c, z_c) = (1600, 1600, 500) \, \mathrm{m}$. 
The initial density and enthalpy are given by  
\eqref{eq:rho_wb} and \eqref{eq:e0}, respectively. The initial velocity field is set to zero everywhere. Impenetrable, free-slip boundary conditions are imposed on all boundaries. The time interval of
interest is $(0, 500]$ s.

We consider a structured uniform mesh of size $h = \Delta x = \Delta y = \Delta z = 40 \, \mathrm{m}$ and set the time step to $\Delta t = 0.1 $ s. Following \cite{marras2013variational}, we set $\mu_a = 12.5$ and $Pr = 1$. 

To generate the ROM model, we sample the $\theta'$ field
every second by collecting a total of 500 snapshots. 
This dataset is split with a T-to-V ratio of $80 \%-20 \%$, i.e.,
we use the first 400 snapshots for training and the remaining 
100 for validation. Given the excellent performance of the
E-CAE network for the density current benchmark, 
we use the same type of network to encode the 3D field 
into $N_d = 4$ latent spaces, i.e., E-CAE network consists of  $L = 4$ convolutional layers with $N_f^1 = 512$, $N_f^2 = 256$, $N_f^3 = 128$, and $N_f^4 = 64$ with
$n_f=3$ sets of filters at each resolution.
The only difference is that for this benchmark we use the 3D CNN layer in the E-CAE network.{ We choose the following hyperparametrs for the RC network:  $N_h = 1200$, $\alpha = 0.015 $ and $\lambda = 0.00055$.} 
%\anna{Do we also use the same RC framework?}

We provide a qualitative comparison of $\theta'$ computed
by the ROM and the FOM in \fig{fig:Vis_RTB3D}, which shows
the results for $t = 340,420,500$ s. Of these three times, 
the first  belongs to the training dataset, while the other 
two belong to the validation set. Just like for the 2D tests, we see that the E-CAE can accurately reconstruct the physical field from its compressed representation
and the combination of E-CAE and RC networks
can accurately predict the spatio-temporal state of the system with remarkable precision. In fact, we observe that
the ROM accurately predicts the system dynamics 
up to the final time (i.e., $t = 500$ s), with only minor differences with respect to the FOM, that are mostly concentrated around the bubble. 

To support the qualitative analysis in \fig{fig:Vis_RTB3D}, \fig{fig:Rec_error_RTB3D} shows the evolution of error \eqref{eq:l2Error} for the reconstruction and prediction phases. 
The error in the reconstruction phase remains around 2.5$\%$, while the prediction error gradually increases to around 8$\%$. 

{We conclude with a comment on the computational times using the machine described at the end of Sec.~\ref{sec:RTB}. 
The FOM simulation for this test requires $4050$ s. The training time of our proposed framework takes  approximately $5400$ s, while the ROM inference time is about $2$ s.}

%\anna{Arash aggiunge un commento sui tempi computazionali.}

\textbf{\begin{figure}[t]
    %\centering \includegraphics[width=100mm,scale=0.6]{Rec_Error_RTB3D}
    \centering \includegraphics[width=.5\textwidth]{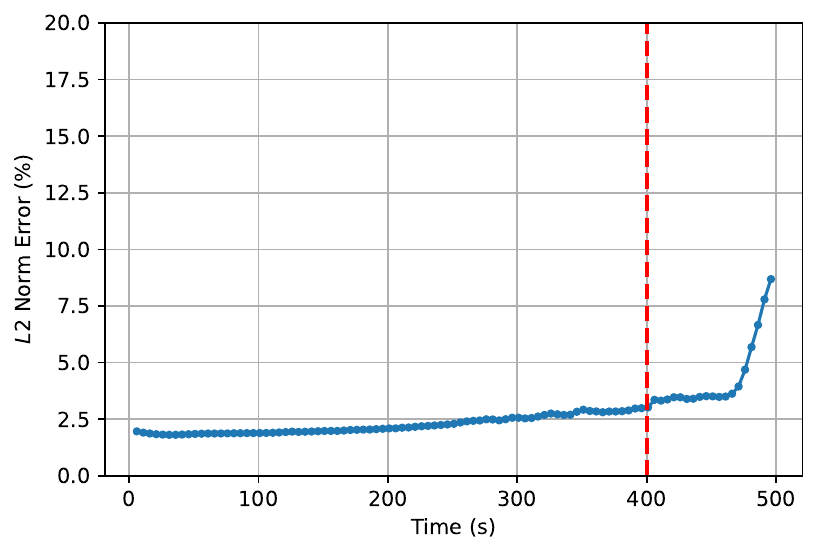}
        \caption{3D rising thermal bubble: evolution of the error \eqref{eq:l2Error}. The red dashed line indicates the beginning of the prediction phase. }   
        \label{fig:Rec_error_RTB3D}
\end{figure}}

\section{Concluding remarks}\label{sec:conc}

Traditional data-driven reduced order models (ROMs), such as Dynamic Mode Decomposition and 
Proper Orthogonal Decomposition with Interpolation, struggle to capture and accurately predict future dynamics of highly nonlinear systems.
We introduced a nonlinear data-driven framework that integrates Convolutional Autoencoder (CAE) and Reservoir Computing (RC) to improve the accuracy of reduced order modeling for mesoscale atmospheric flows. 
Our approach compresses the high-resolution data into a low-dimensional latent space representation through a more capable version of CAE, called Extended Convolutional Autoencoders (E-CAE). Then, a RC network is used to accurately and efficiently predict the future evolution of the latent space.

We assessed our approach with three well-known benchmarks for atmospheric flows, two of which are in two dimensions, while the third is in three dimensions. We show that the E-CAE significantly improves the accuracy of spatial feature extraction when compared to standard CAEs. Moreover, we show that the RC network accurately predicts the future state of the system while significantly reducing the computational cost compared to a full order simulation. This is thanks to the fact that, 
unlike other RNN networks, RC only trains the output weights using a simple linear regression. 
The scalability of our ROM is demonstrated with the 3D benchmark. 

While we believe this study represents a significant step toward a more accurate and efficient approach to accelerate atmospheric flow prediction, improvements are certainly possible. 
One possible improvement may come from the integration of  physics-informed approaches, which are expected to enable accurate predictions over longer time intervals. If one is interested in 
parametric studies, neural operator-based approaches, such as DeepONet \cite{lu2021learning, goswami2022deep} and Fourier Neural Operator \cite{li2020fourier} could be considered to further contain the computational cost.

\section{acknowledgments}\label{sec:akw}

AH and GR acknowledge the support provided by the European Union - NextGenerationEU, in the framework of the iNEST - Interconnected Nord-Est Innovation Ecosystem (iNEST
ECS00000043 – CUP G93C22000610007) consortium and its CC5 Young Researchers initiative. SMACT Industry 4.0 competence center in North-East of Italy is acknowledged by IRISS project initiative.
The authors also like to acknowledge INdAM-GNSC for its support. 

%\printbibliography  
% \bibliographystyle{plain}
%\newpage
\bibliographystyle{unsrt}
\bibliography{mybib}

\end{document}